\documentclass[aos,MSNbibl,nameyear,seceqn,dvips]{arximspdf}

\doi{10.1214/10-AOS865}
\volume{39}
\issue{2}
\pubyear{2011}
\firstpage{1211}
\lastpage{1240}

\makeatletter

\def\bptnote#1{}

\def\N{\mathbb{N}}
\def\R{\mathbb{R}}
\renewcommand{\epsilon}{\varepsilon}

\def\DD{\mathcal D}
\def\UU{\mathcal U}
\def\XX{\mathcal X}
\def\YY{\mathcal Y}

\newtheorem{theorem}{Theorem}[section]
\newtheorem{lem}{Lemma}[section]
\newproclaim{rmk}{Remark}[section]

\renewcommand{\bar}{\overline}
\makeatother

\begin{document}
\begin{frontmatter}

\title{Delta method in large deviations and moderate deviations for estimators}
\runtitle{Delta method in large deviations}

\begin{aug}
\author[A]{\fnms{Fuqing} \snm{Gao}\thanksref{t1}\ead[label=e1]{fqgao@whu.edu.cn}}
\and
\author[B]{\fnms{Xingqiu} \snm{Zhao}\corref{}\thanksref{t2}\ead[label=e2]{xingqiu.zhao@polyu.edu.hk}}
\runauthor{F. Gao and X. Zhao}
\affiliation{Wuhan University
and Hong Kong
Polytechnic University and Zhongnan~University of Economics and Law}
\address[A]{School of Mathematics and Statistics\\
Wuhan University\\
Wuhan 430072\\
China\\
\printead{e1}} 
\address[B]{Department of Applied Mathematics\\
Hong Kong Polytechnic University\\
Hung Hom, Kowloon\\
Hong Kong\\
China\\
and\\
School of Statistics and Mathematics\\
Zhongnan University of Economics\\
\quad and Law\\
Wuhan 430073\\
China\\
\printead{e2}}
\end{aug}
\thankstext{t1}{Supported in part by the
National Natural Science Foundation of China (NSFC) Grant 10871153
and the Specialized Research Fund for the Doctoral Program of Higher
Education (SRFDP) of China Grant 200804860048.}
\thankstext{t2}{Supported in part by the Research Grant Council of
Hong Kong (PolyU 5032/09P).}

\received{\smonth{4} \syear{2010}}
\revised{\smonth{8} \syear{2010}}

%
\begin{abstract}
The delta method is a popular and elementary tool for deriving
limiting distributions of transformed statistics, while applications
of asymptotic distributions do not allow one to obtain desirable
accuracy of approximation for tail probabilities. The large and
moderate deviation theory can achieve this goal. Motivated by the
delta method in weak convergence, a general delta method in large
deviations is proposed. The new method can be widely applied to
driving the moderate deviations of estimators and is illustrated by
examples including the Wilcoxon statistic, the Kaplan--Meier
estimator, the empirical quantile processes and the empirical copula
function. We also improve the existing moderate deviations results
for $M$-estimators and $L$-statistics by the new method. Some
applications of moderate deviations to statistical hypothesis
testing are provided.
\end{abstract}

%
\begin{keyword}[class=AMS]
\kwd[Primary ]{60F10}
\kwd{62G20}
\kwd[; secondary ]{62F12}.
\end{keyword}
\begin{keyword}
\kwd{Delta method}
\kwd{hypothesis testing}
\kwd{Kaplan--Meier estimator}
\kwd{large deviations}
\kwd{$L$-statistics}
\kwd{$M$-estimator}
\kwd{moderate deviations}.
\end{keyword}

\end{frontmatter}

\section{Introduction}\label{sec1}
Consider a family of random variables $\{Y_n,n\geq1\}$ such as the
sample mean. Assume that it satisfies a law of large numbers and
a~fluctuation theorem such as central limit theorem, that is, $Y_n\to
\theta$ in law and there exists a sequence $b_n\to\infty$ such
that $b_n(Y_n-\theta)\rightarrow Y $ in law, where $\theta$ is a
constant and $Y$ is a nontrivial random variable. A large deviation
result~is concerned with estimation of large deviation probabilities
\mbox{$P(|Y_n-\theta|\geq\varepsilon)$} for $\varepsilon>0$. A moderate
deviation result is concerned with estimation of large deviation
probabilities $P(r_n|Y_n-\theta| \geq\varepsilon)$ for
$\varepsilon>0$, where $r_n$ is an intermediate scale between $1$
and $b_n$, that is, $r_n\to\infty$ and $b_n/r_n\to\infty$. In
particular, if $b_n=\sqrt n$, then $r_n=n^{1/2-\delta}$ with
$0<\delta<1$.

The large deviation and moderate deviation problems arise in the
theory of statistical inference quite naturally. For estimation of
unknown parameters and functions, it is first of all important to
minimize the risk of wrong decisions implied by deviations of the
observed values of estimators from the true values of parameters or
functions to be estimated. Such gross errors are precisely the
subject of large deviation theory. The large deviation and moderate
deviation results of estimators can provide us with the rates of
convergence and a useful method for constructing asymptotic
confidence intervals. For the classical large deviation theory with
the empirical measures and sample means, one can refer to \citet
{Sanov}, \citet{GroeneboomOosterhoffRuymgaart} and \citet
{BahadurZabell}. The large deviations for linear combinations of
order statistics ($L$-estimators) were also investigated in
\citet{GroeneboomOosterhoffRuymgaart}. \citet{BahadurZabell}
developed a subadditive method to study the large deviations
for general sample means. For some developments of large deviations
and moderate deviations in statistics, see \citet{Fu-82},
\citet{KK-86}, \citet{SD-87}, \citet{Inglot-Ledwina-90},
\citet{BorovkovMogulskii}, \citet{PS-98}, \citet{Bercu},
\citet{Joutard} and \citet{Arcones-M-ldp-06} for
large deviations of estimators; Kallenberg
(\citeyear{Kallenberg-83}), \citet{gao-spl},
\citet{Arcones-mdp-m-estimator-02}, Inglot and Kallenberg
(\citeyear{Inglot-Kallenberg-03}), \citet{DjelloutGuillinWu} and \citet
{Ermakov-08} for moderate deviations of estimators;
\citet{louani-1}, \citet{worms-2}, \citet{gao},
\citet{lei-wu-05} for large
deviations and moderate deviations of kernel density estimators, and
references therein. On the other hand, large deviations of
estimators can be applied to Bahadur efficiency to determine the
Bahadur slope [\citet{Bahadur-67}, \citet{Nikitin-book},
\citet{HeShao}]
and hypothesis testing [see \citet{Dembo-Z-98}, Sections 3.5
and 7.1].

In statistics, many important estimators are functionals $\Phi(L_n)$
of the
empirical processes $L_n$, and so deriving limiting distribution of $
r_n(\Phi(L_n)-\Phi(\mu))$ from limiting distribution of $
r_n(L_n-\mu)$ is a fundamental problem, where $r_n$ is a sequence
of positive numbers and $\mu$ is the mean of $L_n$. It is well known
that the delta method is
a popular and elementary tool for solving the problem. The method
tells us that the weak convergence of $r_n( X_n-\theta)$ yields the
weak convergence of $r_n(\Phi(X_n)-\Phi(\theta)) $ if $\Phi$ is
Hadamard differentiable (see Section~\ref{sec3}), where $X_n$ is a sequence
of random variables, $\theta$ is a constant and $r_n\to\infty$. For
some developments and applications of the delta method, one can
refer to \citet{Gill-89}, \citet{Kosorok-08}, \citet
{Reeds-76}, and \citet{Vaar-Wellner-96} among others. For example,
\citet{Reeds-76}
systematically developed the use of Hadamard instead of Fr\'{e}chet
differentiability to derive asymptotic distributions of transformed
processes. \citet{AndersenBorganGillKeiding} also described
some applications of the delta method in survival analysis. More
recently, \citet{Vaar-Wellner-96} and \citet{Kosorok-08}
provided an excellent summary of the functional delta method in
terms of a weak convergence.

A natural problem is whether the large deviations of
$r_n(\Phi(X_n)-\Phi(\theta))$ can be obtained from the large
deviations of $ r_n( X_n - \theta)$ if the function~$\Phi$ defined
on a set $\DD_\Phi$ is Hadamard differentiable. When $r_n=r$ for
all $n$~with a constant~$r$, the problem can be solved by the
contraction principle [see \citet{Dembo-Z-98}]. When
$r_n\to\infty$, for each $n\geq1$, define
$\DD_n=\{h;\theta+h/r_n\in\DD_\Phi\}$ and
$f_n(h)=r_n (\Phi ( \theta+h/r_n )-\Phi(\theta
) )$
for all $h\in\DD_n$. Then by Hadamard differentiability, for every
sequence $h_{n}\in\DD_{n}$ converging to~$h$, the\vadjust{\goodbreak} sequence $f_{n}$
satisfies $f_{n}(h_{n})\to\Phi_\theta'(h)$. Note that
$f_n(r_n(X_n-\theta))=r_n(\Phi(X_n)-\Phi(\theta))$. Motivated by
this, we can also consider to use a contraction principle for
establishing the large deviations of $r_n(\Phi(X_n)-\Phi(\theta))$.
However, the existing contraction principles cannot be applicable to
these situations as addressed in Remark~\ref{rem2.1} of next section. For
this reason, we need to extend the contraction principle in large
deviations.

The objective of this paper is to develop a general delta method in
large deviations similar to that in week convergence and applies the
method to solve some moderate deviation problems in statistics. The
remainder of the paper is organized as follows. In Section~\ref{sec2}, we
present an extended contraction principle, while its proof will be
given in the \hyperref[appm]{Appendix}. Then a~general delta method in large deviations
is established by using the extended contraction principle in
Section~\ref{sec3}. In Section~\ref{sec4}, we apply the proposed delta method in large
deviations to some statistical models including censored data,
empirical quantile process, copula function, $M$-estimators and
$L$-statistics. The moderate deviation principles for the Wilcoxon
statistic, the Kaplan--Meier estimator, the empirical quantile
estimator and the empirical copula estimator are established. We
also improve the existing moderate deviation results for
$M$-estimators and $L$-statistics in Section~\ref{sec4}, where our proofs are
different from others but more simple by the new method. Section~\ref{sec5}
presents some applications of the moderate deviation results to
statistical hypothesis testing. Some concluding remarks are made in
Section~\ref{sec6}.

\section{An extended contraction principle}\label{sec2}

As explained in previous section, to establish a delta method in
large deviation, we first need to generalize the contraction
principle in large deviation theory. In this section, we present an extension
of the contraction principle which plays an important role.

First, let us introduce some notation in large deviations [Dembo and~Zei\-touni
(\citeyear{Dembo-Z-98}), \citet{D-stroock-89}]. For a metric space
$\XX$, $\mathcal B(\XX)$ is the Borel $\sigma$-algebra of $\XX$. Let
$(\Omega,\mathcal F, P)$ be a probability space and let~$T$ be an
arbitrary map from $\Omega$ to $\bar{\mathbb R}$, where
$\bar{\mathbb R}=[-\infty,\infty]$ is the space of extended real
numbers. The outer integral of $T$ with respect to $P$ is defined by
\citet{Vaar-Wellner-96}
\[
E^*(T)=\inf\{E(U); U\geq T, U\dvtx \Omega\mapsto\bar{ \mathbb R}\mbox{
measurable and } E(U) \mbox{ exists} \}.
\]
The outer probability of an arbitrary subset $B$ of $\Omega$ is
\[
P^*(B)=\inf\{P(A); A\supset B,A\in\mathcal F\}.
\]
Inner integral and inner probability are defined by
\[
E_*(T)=-E^*(-T) \quad \mbox{and} \quad P_*(B)= 1-P^*(\Omega\setminus B),
\]
respectively.

Let $\{(\Omega_n, \mathcal F_n, P_n),n\geq1\}$ be a sequence of
probability spaces and let $\{X_n, n\geq1\}$ be a sequence of maps
from $\Omega_n$ to $\XX$. Let $\{\lambda(n),n\geq1\}$ be a sequence
of positive numbers tending to $+\infty$ and let $I \dvtx \XX\to[0,
+\infty]$ be inf-compact; that is, $[I\le L]$ is compact for any
$L\in\R$. Then $\{X_n,n\geq1\}$ is said to satisfy the lower bound
of large deviation (LLD) with speed $ \lambda(n) $ and rate
function $I$, if for any open measurable subset $G$ of $\XX$,
%
%
\begin{equation}\label{LLD-def}
l(G):=\liminf_{n\rightarrow\infty}\frac1{\lambda(n)} \log
{P_n}_*(X_n\in G)
\ge- \inf_{x \in G}I(x).
\end{equation}
$\{X_n,n\geq1\}$ is said to satisfy the upper bound of large
deviation (ULD) with speed $\lambda(n)$ and rate function $I$, if
for any closed measurable subset $F$ of $\XX$,
%
%
\begin{equation}\label{ULD-def}
\UU(F):=\limsup_{n\rightarrow\infty} \frac1{\lambda(n)} \log
{P_{n}}^*(X_n\in F) \le-\inf_{x\in F} I(x).
\end{equation}
We say that $\{X_n,n\geq1\}$ satisfies the large deviation
principle (LDP) with speed $\lambda(n)$ and rate function $I$, if
both LLD and ULD hold.

Now, we present the extended contraction principle.

\begin{theorem}[(Extended contraction principle)]\label{ext-contract-princple-thm}
Let $(\XX,d)$ and $(\YY,\rho)$ be two metric spaces. Let
$\{\DD_n,n\geq1\}$ be a sequences of subsets in $(\XX,d)$, and let
$\{f_n\dvtx \DD_n\mapsto\YY; n\in\N\}$ be a family of mappings. Also
for each $n\geq1$, let $ X_n $ be a map from probability space
$(\Omega_n, \mathcal F_n, P_n)$ to $\DD_n$. Suppose that:
\begin{longlist}[(ii)]
\item[(i)] $\{X_n, n\geq1\}$ satisfies the large deviation
principle with speed $\lambda(n)$ and rate function $I$;
\item[(ii)] there exists a mapping
$f\dvtx \{I<\infty\}\mapsto\YY$ such that if for a sequence \{$x_{n}\in
\DD_{n},n\geq1\}$, $x_{n}\to x\in\{I<\infty\}$ as
$n\to\infty$, then
$f_{n}(x_{n})\to f(x)$ as $n\to\infty$.
\end{longlist}
 Then $\{f_n(X_n),n\geq1\}$ satisfies the large deviation
principle with speed $\lambda(n)$ and rate function $I_{ f}$, where
%
%
\begin{equation}\label{ext-contract-princple-thm-eq-3}
I_{ f}(y)=\inf\{I(x) ; f(x)=y\}, \qquad y\in\YY.
\end{equation}
\end{theorem}

The proof of the theorem is given in the \hyperref[appm]{Appendix}.

\begin{rmk}\label{rem2.1} (1) If $\DD_n=\XX$ for all $n\geq1$,
then Theorem~\ref{ext-contract-princple-thm} yields Theorem~2.1 in
\citet{Arcones-LDP-03}. Another popular contraction principle was
given in
Theorem 4.3.23 of \citet{Dembo-Z-98}, in which $\DD_n=\XX$
for all $n\geq1$, $f_n$ is continuous for all $n\geq1$ and for
any $L\in(0,\infty)$,
%
%
\begin{equation}\label{ext-contract-princple-thm-eq-4}
\lim_{n\to\infty}\sup_{x\dvtx I(x)\leq L} \rho(f_n(x),f(x))=0.
\end{equation}
This condition cannot be compared to condition (ii) in Theorem
\ref{ext-contract-princple-thm}.

(2) It is necessary for proving Theorem \ref{Delta-method-thm} to
introduce the sequence of subsets $\DD_n$ in Theorem
\ref{ext-contract-princple-thm}, because subsets $
\{h\in\XX;\theta+{h}/{r_n} \in\DD_\Phi\}, n\geq1$\vadjust{\goodbreak} are not equal,
generally, for $\theta\in\XX$ and a subset $\DD_\Phi$ of a
topological linear spaces $\XX$. In fact, $\DD_\Phi$ is usually a
subset of $\XX$ in applications (see Section~\ref{sec4}).
\end{rmk}

\section{Delta method in large deviations}\label{sec3}

In this section, we establish a delta method in large deviations by
using the extended contraction principle presented in Section~\ref{sec2}.

Let us first recall some conceptions of Hadamard differentiability
[\citet{Gill-89}, \citet{Vaar-Wellner-96}, \citet{Kosorok-08},
\citet{Romisch-05}]. Let $\XX$ and $\YY$ be two metrizable topological
linear spaces. A map $\Phi$ defined on a subset $\DD_\Phi$ of
$\XX$ with values in $\YY$ is called Hadamard differentiable at~$x$
if there exists a continuous mapping $\Phi_x'\dvtx \XX\mapsto\YY$ such
that
%
%
\begin{equation}\label{Hadamard-def-eq-1}
\lim_{n\to\infty}\frac{\Phi(x+t_nh_n)-\Phi(x)}{t_n}=\Phi_x'(h)
\end{equation}
holds for all sequences $t_n$ converging to $0+$ and $h_n$
converging to $h$ in $\XX$ such that $x+ t_nh_n\in\DD_\Phi$ for
every $n$.

\begin{rmk} Linearity of the Hadamard directional derivative
$\Phi_x'(\cdot)$ is not required. In fact, $\Phi_x'(\cdot)$ is often
not linear if $\Phi$ is given by inequality constraints. However, by
the definition, we can see that $\Phi_x'(\cdot)$ is positively
homogenous; that is, $\Phi_x'(t h)=t\Phi_x'( h)$ for all $t\geq0$
and $h\in\XX$.
\end{rmk}

The definition of the Hadamard differentiable may be refined to
Hadamard differentiable tangentially to a set $\DD_0\subset\XX$.
For a subset $\DD_0$ of $\XX$, the map~$\Phi$ is said to be Hadamard
differentiable at $x\in\DD_\Phi$ tangentially to $\DD_0$ if the
limit (\ref{Hadamard-def-eq-1}) exists for all sequences $t_n$
converging to $0+$ and $h_n$ converging to $h$ in $\DD_0$ such that
$x+ t_nh_n\in\DD_\Phi$ for every $n$. In this case, the Hadamard
derivative $\Phi_x'(\cdot)$ is a continuous mapping on $\DD_0$. If
$\DD_0$ is a cone, then $\Phi_x'(\cdot)$ is again positively
homogenous.

\begin{theorem}[(Delta method in large deviation)]\label{Delta-method-thm}
Let $\XX$ and $\YY$ be two metrizable linear topological
spaces and let $d$ and $\rho$ be compatible metrics on $\XX$ and $\YY
$, respectively. Let $\Phi\dvtx \DD_\Phi\subset\XX\mapsto\YY$ be
Hadamard-differentiable at $\theta$ tangentially to $\DD_0$, where
$\DD_\Phi$ and $\DD_0$ are two subsets of $\XX$. Let
$X_n\dvtx \Omega_n\mapsto\DD_\Phi, n\geq1$ be a sequence of maps and
let $r_n,n\geq1 $, be a sequence of positive real numbers
satisfying $r_n\to+\infty$.

If $\{ r_n(X_n-\theta), n\geq1\}$ satisfies the large deviation
principle with speed $\lambda(n)$ and rate function $I$ and
$\{I<\infty\}\subset\DD_0$, then $\{
r_n(\Phi(X_n)-\Phi(\theta)),n\geq1\}$ satisfies the large deviation
principle with speed $\lambda(n)$ and rate function
$I_{\Phi_\theta'}$, where
%
%
\begin{equation}\label{Delta-method-thm-eq-1}
I_{\Phi_\theta'}(y)=\inf\{I(x) ; \Phi_\theta'(x)=y\}, \qquad y\in\YY.
\end{equation}

Furthermore, if $\Phi_\theta'$ is defined and continuous on the
whole space of $\XX$, then $\{
r_n(\Phi(X_n)-\Phi(\theta))-\Phi_\theta'(r_n(X_n-\theta)),n\geq
1\}$
satisfies the large deviation principle with speed $\lambda(n)$ and
rate function
%
%
\begin{equation}\label{Delta-method-thm-eq-2}
I_{\Phi, \theta}(z)=
\cases{\displaystyle 0, &\quad $z=0 $,\cr\displaystyle
+\infty, &\quad otherwise.
}
\end{equation}
In particular, for any $\delta>0$,
%
%
\begin{eqnarray}\label{Delta-method-thm-eq-3}
&&\limsup_{n\to\infty}\frac{1}{\lambda(n)}\log
P_n^* \bigl(\rho\bigl(r_n\bigl(\Phi(X_n)-\Phi(\theta)\bigr)-\Phi_\theta
'\bigl(r_n(X_n-\theta)\bigr),0\bigr)\geq
\delta \bigr)\nonumber
\\[-8pt]
\\[-8pt]
&& \qquad =-\infty.
\nonumber
\end{eqnarray}
\end{theorem}

\begin{pf}
For each $n\geq1$, define
$\DD_n=\{h\in\XX;\theta+h/r_n\in\DD_\Phi\}$ and
\[
f_n\dvtx \DD_n\mapsto\YY, \qquad f_n(h)=r_n\bigl (\Phi (
\theta+h/r_n )-\Phi(\theta) \bigr) \qquad \mbox{for all }
h\in\DD_n.
\]
Then for every sequence $h_{n}\in\DD_{n}$ converging to
$h\in\DD_0$, the sequence $f_{n}$~satis\-fies $f_{n}(h_{n})\to
\Phi_\theta'(h)$. In addition, $\Phi_\theta'(\cdot)$ is continuous
on $\DD_0$. Therefore, Theo\-rem
\ref{ext-contract-princple-thm} implies that
\[
\bigl\{ r_n\bigl(\Phi(X_n)-\Phi(\theta)\bigr),n\geq1\bigr\}
=\bigl\{ f_n\bigl(r_n(X_n-\theta)\bigr),n\geq1\bigr\}
\]
satisfies the large
deviation principle with speed $\lambda(n)$ and rate function
$I_{\Phi_\theta'}$.

Now, we consider the mapping $\varphi_n\dvtx \DD_n\mapsto\YY\times\YY$,
where $\varphi_n(h)=(f_n(h),\break \Phi_\theta'(h))$ for all $h\in\DD_n$.
If $\Phi_{\theta}'(\cdot)$ is continuous on $\XX$, then for every
subsequence $h_{n'}\in\DD_{n'}$ converging to $h\in\XX$,
$\varphi_{n'}(h_{n'})$ converges to $(\Phi_{\theta}'(h),\Phi
_{\theta}'(h))$. Hence, Theorem \ref{ext-contract-princple-thm}
implies $ \{\varphi_n(r_n(X_n-\theta)),n\geq1\}$ satisfies the
large deviation principle with speed $\lambda(n)$ and
rate function
\[
J_{\Phi,\theta}(y_1,y_2)=\inf\{I(x) ; \Phi_\theta'(x)=y_1=y_2\}, \qquad
(y_1,y_2)\in\YY\times\YY.
\]
Therefore, by the classical contraction principle [see \citet
{Dembo-Z-98}, Theorem 4.2.1], we conclude that the difference
\[
\bigl\{r_n\bigl(\Phi(X_n)-\Phi(\theta)\bigr)-\Phi_\theta'\bigl(r_n(X_n-\theta
)\bigr),n\geq1\bigr\}
\]
satisfies the large deviation principle with speed $\lambda(n)$ and
rate function
\[
\inf\{J_{\Phi,\theta}=(y_1,y_2); y_1-y_2=z\}=I_{\Phi, \theta}(z)
 \qquad \mbox{for } z\in\YY.
\]
\upqed
\end{pf}

\section{Moderate deviations of estimators}\label{sec4}

In this section, moderate deviation principles for some estimators
will be established by applying the delta method in large deviation
to Wilcoxon statistic, Kaplan--Meier estimator, the empirical quantile
processes, $M$-estimators
and $L$-statistics.

Let us introduce some notation. Given an arbitrary set $T$ and a
\mbox{Banach}~spa\-ce $(\mathbb B,\|\cdot\|_{\mathbb B})$, the Banach space
$l_\infty(T,\mathbb B)$ is the set of all maps $z\dvtx T\mapsto\mathbb
B$ that are uniformly norm-bounded equipped with the norm
$\|z\|=\sup_{t\in T}\|z(t)\|_{\mathbb B}$. Let $l_{\infty}(T)$ be
the Banach space of all bounded real functions $x$ on $T$, equipped
with the sup-norm $\|x\|=\sup_{t\in T}|x(t)|$. It is a nonseparable
Banach space if $T$ is infinite. On $l_{\infty}(T)$, we will
consider the $\sigma$-field $\mathcal B$ generated by all balls and
all coordinates $x(t), t\in T$.

Let $(S,d)$ be a complete separable and measurable metric space and
let $b\mathcal S$ be the space of all bounded real measurable
functions on $(S, \mathcal S)$ where $\mathcal S$ is the Borel
$\sigma$-algebra of $S$. Let $\{X, X_n,n\geq1\}$ be a sequence of
i.i.d. random variables with values in $S$ on a probability space
$(\Omega,\mathcal F,P)$, of law $\mu$. Let $L_n$ denote the
empirical measures; that is,
\[
L_n=\frac{1}{n}\sum_{i=1}^n \delta_{X_i}, \qquad n\geq1.
\]

For given a class of functions $\mathfrak{F}\subset b\mathcal S$,
let $l_\infty(\mathfrak{F})$ be the space of all bounded real
functions on $\mathfrak{F}$ with sup-norm
$\|F\|_{\mathfrak{F}}=\sup_{f\in\mathfrak{F}}|F(f)|$. This is a
Banach space. Every $\nu\in M_b(S)$ [the space of signed measures of
finite variation on $(S,\mathcal S)$] corresponds to an element
$\nu^{\mathfrak{F}}=\nu(f)=\int f \,d\nu$ for all $f\in\mathfrak{F}$.

Let $D[a,b]$ denote the Banach space of all right continuous with
left-hand limits functions $z\dvtx [a,b]\mapsto\mathbb R$ on an
interval $[a,b]\subset\overline{\mathbb R}$ equipped with the
uniform norm. Let $BV[a,b]$ denote the set of all cadlag functions
with finite total variation and set $ BV_M[a,b]= \{A\in
BV[a,b]; \int|dA|\leq M \}, $ where the notation $\int|dA|$
denotes the total variation of the function $A$. In this article, we
also let $\{a_n=a(n),n\geq1\}$ be a sequence of real numbers such
that as $n\to\infty$,
\[
a_n \to\infty\quad\mbox{and} \quad a_n/\sqrt{n}\to0.
\]

\subsection{Moderate deviations for Wilcoxon statistic}\label{sec4.1}

 Let $X_1, \ldots, X_m$ and $Y_1,\break\ldots, Y_n$ be independent
samples from distribution functions $F$ and $G$ on $\mathbb R$,
respectively. If $F_m$ and $G_n$ are the empirical distribution
functions of the two samples; that is,
\[
F_m(x)=\frac{1}{m}\sum_{i=1}^m \delta_{X_i}((-\infty,x])
 \quad \mbox{and} \quad G_n(x)=\frac{1}{n}\sum_{i=1}^n
\delta_{Y_i}((-\infty,x]),
\]
then the Wilcoxon statistic is defined by
$
W_{m,n}=\int F_m \,dG_n.
$
It is an estimator of $P(X\leq Y)$.

\begin{theorem}\label{wilcoxon-mdp-thm}
Assume that $m/(m+n)\to\lambda\in(0,1)$ as $m,n\to\infty$. Then
%
%
\begin{equation}\label{wilcoxon-mdp-thm-eq-1}
\biggl \{ \frac{\sqrt{mn/(m+n)}}{a(mn/(m+n))} \biggl(\int F_m \,dG_n-\int F\,dG \biggr),n\geq1 \biggr\}
\end{equation}
satisfies the LDP in $ \mathbb R$ with speed $a^2(mn/(m+n))$ and
rate function $I^W$ defined by
%
%
\begin{equation}\label{wilcoxon-mdp-thm-eq-rate}
I^{W}(x)=\frac{x^2}{2 (\lambda \operatorname{Var}(F(Y))+(1-\lambda)\operatorname{Var}(G(X)))}.
\end{equation}
\end{theorem}

\begin{pf}
Applying Theorem 2 of \citet{Wu-94} to
$L_n^X=\frac{1}{n}\sum_{i=1}^n \delta_{X_i}$,
$\mathfrak{F}_1=\{(-\infty,x];x\in\mathbb R\}$, and $L_n^Y=
\frac{1}{n}\sum_{i=1}^n \delta_{Y_i}$,
$\mathfrak{F}_2=\{(-\infty,y];y\in\mathbb R\}$, respectively, and
using the product principle in large deviations [\citet
{Dembo-Z-98}], we obtain that $ \{ \frac{\sqrt{n}}{a_n}
 (F_n-F,G_n-G ),n\geq1 \} $ satisfies the LDP in
$l_\infty(\mathbb R ) \times l_\infty(\mathbb R )$ with speed
$a_n^2$ and rate function $ \{I_F(\alpha)+ I_G(\beta)\}, $ where
\begin{eqnarray*}
I_F(\alpha)
&=&\inf
\biggl \{
\frac{1}{2}
\int\gamma^2(x)\,dF(x);
 \ \int\gamma(x)\,dF(x)=0,\alpha(t)=\int_{(-\infty,t]}\gamma(x)\,dF(x)\\
 &&\hphantom{\inf
\biggl \{}
 \hspace*{108pt}\mbox{for each }t\in\mathbb R, \ \gamma\dvtx \mathbb R\to\mathbb
R \mbox{ is measurable}
\!\biggr\}\\
&=&
\cases{\displaystyle \frac{1}{2}\int|\alpha'
_F(x)|^2\,dF(x), &\quad if $\alpha\ll F$ and
$\displaystyle\lim_{|t|\to\infty}|\alpha(t)|=0$, \vspace*{2pt}\cr\displaystyle
\infty, &\quad otherwise
}
\end{eqnarray*}
and $\alpha_F'=d\alpha/dF$. Since
$\frac{\sqrt{m}}{\sqrt{mn/(m+n)}}\to(1-\lambda)^{1/2}$ and $ \frac
{\sqrt{n}}{\sqrt{mn/(m+n)}}\to
\lambda^{1/2},$ then
\[
\biggl \{\frac{\sqrt{mn/(m+n)}}{a(mn/(m+n))} (F_m-F,G_n-G
),n\geq1 \biggr\}
\]
satisfies the LDP in $l_\infty(\mathbb R ) \times l_\infty(\mathbb R
)$ with speed $ {a^2(mn/(m+n))} $ and rate function given by
\[
I_{F,G}(\alpha,\beta)=\frac{1}{1-\lambda}I_F(\alpha)+\frac{1}{
\lambda}I_G(\beta).
\]

Note that $\{I_{F,G}(\alpha,\beta)<\infty\}\subset BV (\mathbb
R)\times BV (\mathbb R)$ and $(F_m,G_n )\in BV_1(\mathbb R)\times
BV_1(\mathbb R).$ For each $M\geq1$, we consider the map
$\Phi\dvtx D(\mathbb R)\times BV_M(\mathbb R)\mapsto\mathbb R$ defined
as
\[
\Phi(A,B)=\int_{\mathbb R} A(s)\,dB(s).
\]
Then $ \Phi(F_m,G_n)=\int F_m \,dG_n$, and by Lemma 3.9.17 of \citet
{Vaar-Wellner-96}, $\Phi$ is Hadamard differentiable at each
$(A,B)\in\DD_\Phi=\{\int|dA|<\infty\}$ and the derivative is given
by
\[
\Phi'_{A,B}(\alpha,\beta)=\int_{\mathbb R} A(s)\,d\beta(s)+\int_{\mathbb R} \alpha(s)\,dB(s),
\]
where $\int_{(a,b]} A(s)\,d\beta(s) $ is defined via integration by
parts if $\beta$ is not of bounded variation; that is,
\[
\int_{(a,b]} A(s)\,d\beta(s) =A(b)\beta(b)-A(a)\beta(a)-\int_{(a,b]}
\beta(s-) A(s).\
\]
Thus, by Theorem \ref{Delta-method-thm} with
$\DD_0=\{(\alpha,\beta); I_F(\alpha)<\infty,I_G(\beta)<\infty\}
$, we
conclude that
\[
 \biggl\{\frac{\sqrt{mn/(m+n)}}{a(mn/(m+n))} \biggl(\int F_m \,dG_n-\int F
\,dG \biggr),n\geq1 \biggr\}
\]
satisfies the LDP on $\mathbb R$ with
speed $ {a^2(mn/(m+n))} $ and rate function given~by
\begin{eqnarray*}
I^W(x)&=&\inf \biggl\{\frac{1}{1-\lambda}I_F(\alpha)+\frac{1}{
\lambda
}I_G(\beta),\int F(s)\,d\beta(s)+\int\alpha(s)\,dG(s)=x \biggr\}\\
&=&\inf \biggl\{\frac{1}{2(1-\lambda)}\int(\alpha'_F)^2\,dF+\frac{1}{2
\lambda} \int(\beta'_G)^2\,dG,
\\
&&\hphantom{\inf \biggl\{} \int F \beta_G'\,dG-\int
G\alpha_F'
\,dF =x, \alpha\ll F,
\\
&&\hphantom{\inf \biggl\{}\beta\ll G, \lim_{|t|\to\infty}|\alpha(t)|=0, \lim_{|t|\to\infty
}|\beta(t)|=0
 \biggr\}\\
&=&\frac{x^2}{2 (\lambda \operatorname{Var}(F(Y))+(1-\lambda)\operatorname{Var}(G(X)))}.
\end{eqnarray*}
\upqed
\end{pf}

\subsection{Moderate deviations for Kaplan--Meier estimator}\label{sec4.2}

Let $X$ and $C$ be independent, nonnegative random variables with
distribution functions~$F$ and $G$. Let $X_1, \ldots, X_n$
be i.i.d. random variables distributed according to the
distribution function~$F$ and let $C_1,\ldots, C_n$ be i.i.d. random
variables distributed according to the distribution function $G$.
$X_1, \ldots, X_n$ and $C_1, \ldots, C_n$ are assumed to be
independent. Observed data are the pairs $(Z_1,\Delta_1), \ldots,
(Z_n,\break\Delta_n),$ where $Z_i=X_i\wedge C_i$, and
$\Delta_i=1_{\{X_i\leq C_i\}}$. The cumulative hazard function is
defined by
%
%
\begin{equation}\label{cumulative-hazard-function}
\Lambda(t)=\int_{[0,t]}
\frac{1}{\overline{F}(s)}\,dF(s)=\int_{[0,t]}\frac{1}{\overline
{H}(s)}\,dH^{uc}(s),
\end{equation}
where
$
\overline{F}(t)=P(X\geq t) \mbox{ and } \overline{H}(t)=P(Z\geq
t)
$
are (left-continuous) survival distributions,
and $H^{uc}(t)=P(Z\leq
t,\Delta=1)$ is a subdistribution function of the uncensored
observations, where $\Delta=1_{\{X\leq C\}}$. We also denote
$H^{c}(t)=P(Z\leq t, \Delta=0)$. The Nelson--Aalen estimator is
defined by
%
%
\begin{equation}\label{Nelson-Aalen-estimator}
\Lambda_n(t)=\int_{[0,t]} \frac{1}{\overline{H}_n(s)}\,dH_n^{uc}(s),
\end{equation}
where
%
%
\begin{equation}\label{uncensored-survival-distributions}
H_n^{uc}(t)=\frac{1}{n}\sum_{i=1}^n 1_{\{Z_i\leq t,\Delta_i=1\}}
\quad\mbox{and}\quad\overline{H}_n(t)=\frac{1}{n}\sum_{i=1}^n
1_{\{Z_i\geq t\}}
\end{equation}
are the empirical subdistribution functions of the uncensored
failure time and the survival function of the observation times,
respectively.

The distribution function $F(t)$ can be rewritten as
\[
1-F(t)=\prod_{0<s\leq t} \bigl(1-d\Lambda(s) \bigr).
\]
The Kaplan--Meier estimator $\hat{F}_n(t)$ for the distribution
function $F(t)$ is defined by
%
%
\begin{equation}\label{Kaplan-Meier-estimator}
1-\hat{F}_n(t)=\prod_{0<s\leq t} \bigl(1-d\Lambda_n(s) \bigr).
\end{equation}
The Kaplan--Meier estimator $\hat{F}_n$ is the nonparametric maximum
likelihood estimator of $F$ in the right censored data model,
proposed by \citet{Kaplan-Meier-58}. \citet{Dinwoodie-93}
studied large
deviations for censored data and established a large deviation
principle for $\sup_{x\in\tau}|\hat{F}_n(x)-F(x)|$ where $\tau$ is
a fixed time satisfying $\{1-F(\tau)\}\{1-G(\tau)\}>0$. \citet
{BLM-99} obtained an exponential inequality for
$\sup_{x\in\mathbb R}\{(1-G(x))|\hat{F}_n(x)-F(x)|\}$. \citet
{Wellner-07} provided a bound for the constant in the inequality. In this
subsection, we establish its moderate deviation principle.

\begin{theorem}\label{Nelson-Aalen-mdp-thm} Let $\tau>0$ such that
$H(\tau)<1$. Then \mbox{$\{\frac{\sqrt{n}}{a(n)} (\Lambda_n-\Lambda),n\geq1 \}$}
satisfies the LDP in $ D[0,\tau]$ with speed $a^2(n)$
and rate function $I^\Lambda$ given by
%
%
\begin{eqnarray}\label{Nelson-Aalen-mdp-thm-rate-1}
 I^\Lambda(\phi)&=&\inf
\biggl\{
I_{F,G}(\alpha,\beta);
\
\int_{[0,t]}
\frac{1}{\overline{H}(s)}\,d\alpha(s)-\int_{[0,t]}
\frac{\beta(s)}{\overline{H}^2(s)} \,d {H}^{uc}(s) =\phi(t),\nonumber
\\[-8pt]
\\[-8pt]
&&\hspace*{198.5pt}\hphantom{\inf
\biggl\{}\mbox{for
any } t\in[0,\tau]
\biggr\},
\nonumber
\end{eqnarray}
where
%
%
\begin{equation}\label{Nelson-Aalen-mdp-thm-rate-2}
 I_{F,G}(\alpha,\beta)=
\cases{\displaystyle
\frac{1}{2} \biggl(\int
|\alpha'_{H^{uc}}(u)|^2\,dH^{uc}(u)+\int|{(\alpha+\beta)'}
_{H^{c}}(u)|^2\,dH^{c}(u) \!\biggr)\!,\vspace*{1pt}\cr\displaystyle
 \hspace*{16pt}\qquad \mbox{if } \alpha\ll
H^{uc},\ \alpha+\beta\ll H^c \mbox{ and }
\lim_{t\to\infty}|\beta(t)|=0,\vspace*{2pt}\cr\displaystyle
\infty, \qquad \mbox{otherwise.}}\hspace*{-30pt}
\vspace{-5pt}
\end{equation}
\end{theorem}

\begin{pf}
The pair $(H_n^{uc},\overline{H}_n)$ can be identified with the
empirical distribution of the observations indexed by the functions
$ \mathfrak{F}_1=\{I_{\{z\leq t,\Delta=1\}},t\in\mathbb R\}$ and
$\mathfrak{F}_2=\{I_{\{z\geq t\}}, t\in\mathbb R\}. $ It is easy to
verify that the two classes $\mathfrak{F}_1$ and $\mathfrak{F}_2$
are Donsker classes and the mapping $\Psi\dvtx l_\infty(\mathfrak{F} )
\mapsto l_\infty(\mathfrak{F}_1)\times l_\infty(\mathfrak{F}_2) $
defined by $ \phi\longrightarrow(\phi|_{\mathfrak{F}_1},
\phi|_{\mathfrak{F}_2}) $ is continuous, where
$\mathfrak{F}=\bigcup_{j=1}^2 \mathfrak{F}_j$. Applying Theorem~2 of
\citet{Wu-94} to $L_n=\frac{1}{n}\sum_{i=1}^n\delta_{(Z_i,\Delta
_i)}$ and
$\mathfrak{F}$, and the classical contraction principle [see \citet
{Dembo-Z-98}, Theorem 4.2.1] to $ \Psi$, we can get that\looseness=-1
\[
 \biggl\{\frac{\sqrt{n}}{a(n)} (H_n^{uc}-H^{uc},
\overline{H}_n-\overline{H} ),n\geq1 \biggr\}
\]
satisfies the LDP on $D([0,\tau])\times D([0,\tau])$ with speed $
{a^2(n)} $ and rate function
\begin{eqnarray*}
I_{F,G}(\alpha,\beta)&=&\inf \biggl\{
 \frac{1}{2} \biggl(\int
\gamma_1^2(u)\,dH^{uc}(u)+\int\gamma_0^2(u)\,dH^{c}(u) \biggr);\\
&&\hphantom{\inf \biggl\{}
\int\gamma_1(u)\,dH^{uc}(u) +\int\gamma_0(u)\,dH^c(u)=0,\\
&&\hphantom{\inf \biggl\{}
\mbox{and for any } t\in[0,\infty),
\int_{[0,t]}\gamma_1(u)\,dH^{uc}(u)=\alpha(t),\\
&&\hphantom{\inf \biggl\{}
\int_{[t,\infty)} \gamma_1(u)\,dH^{uc}(u)+\int_{[t,\infty)}
\gamma_0(u)\,dH^{c}(u)=\beta(t)
 \biggr\}\\
&=&
\cases{\displaystyle
\frac{1}{2} \biggl(\int
|\alpha'_{H^{uc}}(u)|^2\,dH^{uc}(u)+\int|{(\alpha+\beta)'}
_{H^{c}}(u)|^2\,dH^{c}(u) \biggr),\vspace*{1pt}\cr\displaystyle
 \hspace*{16pt}\qquad \mbox{if } \alpha\ll
H^{uc},\ \alpha+\beta\ll H^c \mbox{ and }
\lim_{t\to\infty}|\beta(t)|=0,\vspace*{2pt}\cr
\infty, \qquad \mbox{otherwise.}
}
\end{eqnarray*}
Set
$
\DD_\Phi= \{(A,B)\in BV _1([0,\tau])\times
D([0,\tau]); B\geq\overline{H}(\tau)/2\}.
$
By the Dvoretzky--Kiefer--Wolfowitz inequality [cf. \citet
{Massart-90}], 
for any
$\epsilon>0$,
\[
P \Bigl( \sup_{t\in
[0,\tau]}|\overline{H}_n(t)-\overline{H}(t)|>\epsilon \Bigr)\leq2
\exp\{- 2 n\epsilon^2\}.
\]
In particular, take $\epsilon=\overline{H}(\tau)/2$, then we have
%
%
\begin{eqnarray}\label{Kaplan-Meier-mdp-exp-eq-eq}
&&\limsup_{n\to\infty} \frac{1}{a^2(n)} \log P^* \bigl(
(H_n^{uc},\overline{H}_n)\notin\DD_\Phi \bigr) \nonumber
\\[-8pt]
\\[-8pt]
&& \qquad \leq
\limsup_{n\to\infty} \frac{1}{a^2(n)} \log P^* \bigl(
\overline{H}_n(t)\leq\overline{H(\tau)}/2 \bigr)=-\infty.
\nonumber
\end{eqnarray}

Consider the maps $\Phi_1\!\dvtx\!\DD_\Phi\!\subset\!
BV_1([0,\tau])\!\times\!
D([0,\tau])\!\mapsto\!BV_1([0,\tau])\!\times\! D([0,\tau])$ and
$\Phi_2\!\dvtx\!
BV([0,\tau])\!\times\! D([0,\tau])\!\mapsto\! D([0,\tau])$
defined as
\[
\Phi_1(A,B)=(A,1/B) \quad\mbox{and}\quad
\Phi_2\dvtx (A,B)\mapsto\int_{[0,\cdot]} B\,dA.
\]
Define $\Phi(A,B)=\Phi_2(\Phi_1(A,B))$. Then $
\Phi(H_n^{uc},\overline{H}_n)=\Lambda_n$, $ \Phi(H^{uc},\overline
{H})=\Lambda$
and by Lemma 3.9.17 of \citet{Vaar-Wellner-96}, $\Phi$ is
Hadamard differentiable at each $(A,B)\in\DD_\Phi$. The derivative
is given by
\[
\Phi'_{A,B}(\alpha,\beta)(t)=\int_{[0,t]} \frac{1}{B(s)}\,d\alpha(s)-\int_{[0,t]} \frac{\beta(s)}{B^2(s)}\,dA(s).
\]
Applying Theorem \ref{Delta-method-thm} to
$\Omega_n=\{(H_n^{uc},\overline{H}_n) \in\DD_\Phi\}$,
$P_n(\cdot)=P( \cdot|\Omega_n)$ and $\DD_0=\DD_\Phi$ together with
(\ref{Kaplan-Meier-mdp-exp-eq-eq}), we conclude that $
\{\frac{\sqrt n}{a(n)} (\Lambda_n-\Lambda),n\geq1 \}
$
satisfies the LDP in $
D[0,\tau]$ with speed $a^2(n)$ and rate function $I^\Lambda$ given
by
\begin{eqnarray*}
&&\hspace*{-4pt} I^\Lambda(\phi)=\inf
 \biggl\{I_{F,G}(\alpha,\beta);\
 \int_{[0,t]} \frac{1}{\overline{H}(s)}\,d\alpha
(s)-\int_{[0,t]}
\frac{\beta(s)}{\overline{H}^2(s)} \,d
{H}^{uc}(s) =\phi(t), \qquad \qquad \\
&&\hspace*{-4pt} \hspace*{198pt}\hphantom{I^\Lambda(\phi)=\inf
 \biggl\{} \mbox{for any } t\in
[0,\tau]
%
 \biggr\}. \qquad \qquad 
\end{eqnarray*}
\upqed
\end{pf}

Next, we give some other representations. Let
$\{(G^{uc}(t),\overline{G}(t)),t\in[0,\tau]\}$ be a zero-mean
Gaussian process with covariance structure
\begin{eqnarray*}
E(G^{uc}(s)G^{uc}(t))&=&H^{uc}(s\wedge t )-H^{uc}(s)H^{uc}(t),\\
E(\overline{G}(s)\overline{G}(t ))&=&\overline{H}(s\vee t )
-\overline{H}(s)\overline{H}(t),
\end{eqnarray*}
 and
 \[
E(G^{uc}(s)\overline{G}(t) )=\bigl(H^{uc}(s )-H^{uc}(t-
)\bigr)I_{(-\infty,s]}(t)-H^{uc}(s)\overline{H}(t).
\]
Set $\tilde{T}=\{(j,t),j=1,2, t\in[0,\tau]\}$ and
\[
\tilde{Z}=\bigl\{\tilde{Z}_{(j,t)}; j=1,2, t\in[0,\tau],
\tilde{Z}_{(1,t)}=G^{uc}(t),\tilde{Z}_{(2,t)}=\overline{G}(t)\bigr\}.
\]
Then by Theorem 5.2 of \citet{Arcones-LDP-04}, $
\{\{\tilde{Z}_{(j,t)}/\sqrt{\lambda(n)}, (j,t)\in\tilde{T}\}
,n\geq
1\}$ satisfies LDP on $l_{\infty}(\tilde{T})$ with speed
$\lambda(n)$ and rate function given by
\[
\tilde{I}(x)=\inf
 \bigl\{ \tfrac{1}{2} E (\gamma^2 );\
 \gamma\in\mathcal L,\ E \bigl(
\gamma\tilde{Z}_{(j,t)} \bigr)=x_{(j,t)} \mbox{ for all } (j,t)\in
\tilde{T}
 \bigr\},
\]
where $\mathcal L$ is the closed vector space of $L^2(P)$ generated
by $\{\tilde{Z}_{(j,t)},(j,t)\in\tilde{T}\}$. Since the mapping
$\Psi\dvtx
l_\infty(\tilde{T}) \mapsto l_\infty([0,\tau],\mathbb R^2) $
defined by
\[
\bigl\{\phi_{(j,t)},(j,t)\in\tilde{T}\bigr\}\longrightarrow
\bigl\{\bigl(\phi_{(1,t)},\phi_{(2,t)}\bigr), t\in[0,\tau]\bigr\}
\]
is continuous, then by the classical contraction principle [see
\citet{Dembo-Z-98}, Theorem 4.2.1], we know that
$ \{\frac{1}{\sqrt{\lambda(n)}}(G^{uc},\overline{G}), n\geq
1 \}$ satisfies the LDP on $D([0,\tau])\times D([0,\tau])$ with
speed $\lambda(n)$ and rate function $I_{F,G}(\alpha,\beta)$, where
$\lambda(n)\to\infty$ as $n\to\infty$.\vspace{1pt}

Define
$
M^{uc}(t)=G^{uc}(t)-\int_{[0,t]} \overline{G}(u)\,d\Lambda(u)
$
and
%
%
\begin{equation}\label{Kaplan-Meier-mdp-Z-def}
Z(t)=\int_{[0,t]}\frac{1}{\overline{H}(s)}\,dG^{uc}(s)-\int_{[0,t]}
\frac{\overline{G}(s)}{\overline{H}^2(s)} \,d {H}^{uc}(s),
\end{equation}
where the first term on the right-hand side is to be understood via
integration by parts. Then $M^{uc}$ is a zero-mean Gaussian
martingale with covariance function [\citet{Vaar-Wellner-96}, page 384]
\[
E(M^{uc}(s)M^{uc}(t))=\int_{[0,s\wedge t]} \overline{H}(u)\bigl(1-\Delta
\Lambda(u)\bigr)\,d\Lambda(u),
\]
where $\Delta\Lambda(u)=\Lambda(u)-\Lambda(u-) $ and
$
Z(t)=\int_{[0,t]}\frac{1}{\overline{H}(s)}\,dM^{uc}(s)
$
is a
zero-mean Gaussian process with covariance function
\[
E(Z(s)Z(t))=\int_{[0,s\wedge t]}\frac{1-\Delta\Lambda(u)}{
\overline{H}(u)}\,d\Lambda(u).
\]
Therefore, by Theorem \ref{ext-contract-princple-thm}, we conclude
that \mbox{$\{\{Z(t)/\sqrt{\lambda(n)},t\in[0,\tau]\}, n\geq1\}$}
satisfies the LDP on
$D([0,\tau]) $ with speed $\lambda(n)$ and rate function
$I^{\Lambda}(\phi)$. Furthermore, from Theorem 5.2 of \citet
{Arcones-LDP-04}
and Theorem 3.1 of \citet{Arcones-LDP-03},
we have the following result.

\begin{theorem}\label{Nelson-Aalen-mdp-thm-2} Let $\tau>0$ such that
$H(\tau)<1$. Then $
\{\frac{\sqrt{n}}{a(n)} (\Lambda_n-\Lambda ),n\!\geq\!1 \}
$
satisfies the LDP in $ D[0,\tau]$ with speed
$a^2(n)$ and rate function $I^{\Lambda}$ given by
%
%
\begin{eqnarray}\label{Nelson-Aalen-mdp-thm-rate-2-1}
 I^\Lambda(\phi)&=&\sup_{m\geq1,t_1,\ldots, t_m\in
[0,\tau]}\sup_{\alpha_1,\ldots, \alpha_m\in\mathbb R} \Biggl\{
\sum_{i=1}^m \phi_{t_i}\alpha_i \nonumber\\
&&\hphantom{\sup_{m\geq1,t_1,\ldots, t_m\in
[0,\tau]}\sup_{\alpha_1,\ldots, \alpha_m\in\mathbb R} \Biggl\{}
{}-\frac{1}{2} \sum_{k,j=1}^m
\alpha_{k}\alpha_{j}
\\
&&\hphantom{\sup_{m\geq1,t_1,\ldots, t_m\in
[0,\tau]}\sup_{\alpha_1,\ldots, \alpha_m\in\mathbb R} \Biggl\{
{}-\frac{1}{2} \sum_{k,j=1}^m}
{}\times \int_{[0,t_k\wedge t_j]}\frac{1-\Delta
\Lambda(u)}{ \overline{H}(u)}\,d\Lambda(u) \Biggr\}.
\nonumber
\end{eqnarray}
In particular, for any $r>0$,
%
%
\begin{equation}\label{Nelson-Aalen-mdp-thm-rate-2-2}
\lim_{n\to\infty}\frac{1}{a^2(n)}\log
P \biggl(\frac{\sqrt{n}}{a(n)}\sup_{x\in[0,\tau]}|\Lambda
_n(x)-\Lambda(x)|\geq
r \biggr)=- \frac{r^2}{2\sigma_\Lambda^2},
\end{equation}
where
$
\sigma_\Lambda^2=\int_{[0,\tau]}\frac{1-\Delta\Lambda(u)}{
\overline{H}(u)}\,d\Lambda(u).
$
\end{theorem}

Now we present the moderate deviations for the Kaplan--Meier
estimator $\hat{F}_n(t)$.

\begin{theorem}\label{Kaplan-Meier-mdp-thm} Let $\tau>0$ such that
$H(\tau)<1$. Then
\mbox{$\{\frac{\sqrt{n}}{a(n)} (\hat{F}_n-F ),n\geq1 \}$}
satisfies the LDP in $ D[0,\tau]$ with speed $a^2(n)$ and rate
function $I^{KM}$ given~by
%
%
\begin{eqnarray}\label{Kaplan-Meier-mdp-thm-rate-f}
&&\hspace*{-5pt} I^{KM}(\phi) = \sup_{m\geq 1,t_1,\ldots,t_m\in
[0,\tau]}\sup_{\alpha_1, \ldots, \alpha_m\in\mathbb R} \Biggl\{
\sum_{i=1}^m \phi_{t_i}\alpha_i \nonumber\\
&&\hspace*{-5pt} \hphantom{I^{KM}(\phi) =\sup_{m\geq 1,t_1,\ldots,t_m\in
[0,\tau]}\sup_{\alpha_1, \ldots, \alpha_m\in\mathbb R} \Biggl\{}
{}-\frac{1}{2} \sum_{k,j=1}^m
\alpha_{k}\alpha_{j}\nonumber
\\[-8pt]
\\[-8pt]
&&\hspace*{157pt}
{}\times\int_{[0,t_k\wedge
t_j]}\frac{(1-F(t_k))(1-F(t_j))}{(1-\Delta
\Lambda(u))\overline{H}(u)}\,d\Lambda(u) \Biggr\}.
\nonumber
\end{eqnarray}
%
In particular, for any $r>0$,
%
%
\begin{equation}\label{Kaplan-Meier-mdp-thm-cor-eq-1}
\lim_{n\to\infty}\frac{1}{a^2(n)}\log
P\biggl (\frac{\sqrt{n}}{a(n)}\sup_{x\in[0,\tau]}|\hat
{F}_n(x)-F(x)|\geq
r \biggr)=- \frac{r^2}{2\sigma_{KM}^2},
\end{equation}
where
\[
\sigma_{KM}^2=\sup_{t\in[0,\tau]}\bigl(1-F(t)\bigr)^2\int_{(0,
t]}\frac{1}{(1-\Delta\Lambda(u))\overline{H}(u)}\,d\Lambda(u).
\]
\end{theorem}

\begin{pf} The map $\Phi\dvtx BV[0,\tau]\subset D[0,\tau]\mapsto
D[0,\tau]$ is defined as
\[
\Phi(A)(t)=\prod_{0<s\leq t} \bigl(1+dA(s)\bigr).
\]
Then,
$
1-F(x)=\Phi(-\Lambda)(x)
 \mbox{ and } 1-\hat{F}_n(x)=\Phi(-\Lambda_n)(x).
$
Since $H(\tau)<1$, there exists some $M\in(0,\infty)$ such that
$\Lambda\in BV_M[0,\tau]$. From
(\ref{Nelson-Aalen-mdp-thm-rate-2-2}), we have
\begin{eqnarray*}
&&\limsup_{n\to\infty}\frac{1}{a^2(n)}\log P^* (
\Lambda_n \notin BV_{M+1}[0,\tau] )\\
&& \qquad \leq \lim_{n\to\infty}\frac{1}{a^2(n)} \log P \Bigl(
\sup_{x\in[0,\tau]}|\Lambda_n(x)-\Lambda(x)|\geq1 \Bigr)=-
\infty.
\end{eqnarray*}
By Lemma 3.9.30 of \citet{Vaar-Wellner-96}, we know that
$\Phi$ is Hadamard differentiable in $ BV_{M+1}[0,\tau]$ with
derivative
\[
\Phi'_{A}(\alpha)(t)=\int_{(0,t]}
\Phi(A)(0,u)\Phi(A)(u,t]\,d\alpha(u),
\]
where
$
\Phi(A)(u,t]=\prod_{u<s\leq t} (1+dA(s)).
$
Applying Theorem \ref{Delta-method-thm} to\vspace*{1pt} $\Omega_n=\{\Lambda_n \in
BV_{M+1}[0,\tau] \}$, $P_n(\cdot)=P( \cdot|\Omega_n)$ and $\DD_0=
BV_{M+1}[0,\tau]$, we obtain from Theo\-rem~\ref{Nelson-Aalen-mdp-thm} that $ \{\frac{\sqrt{n}}{a(n)}
(\hat{F}_n-F ),n\geq1 \} $ satisfies the LDP in $ D[0,\tau]$ with
speed $a^2(n)$ and rate function $ \tilde{I}^{KM}$ given by
\begin{eqnarray*}
\tilde{I}^{KM}&=&\inf \biggl\{
I^{\Lambda}(\alpha);\\
&&\hphantom{\inf \biggl\{} \int_{(0,t]}
\Phi(F)(0,u)\Phi(F)(u,t]\,d\alpha(u)=\phi(t), \mbox{ for any } t\in
[0,\tau]
 \biggr\}.
\end{eqnarray*}

On the other hand, we consider the process
$\Phi'_{-\Lambda}(Z)(t)$, where $Z$ is defined by
(\ref{Kaplan-Meier-mdp-Z-def}). Since
\begin{eqnarray*}
\Phi'_{-\Lambda}(Z)(t)&=&\int_{(0,t]}\frac{(1-F(u-))
((1-F(t))}{1-F(u)}\,dZ(u)\\
&=&\bigl(1-F(t)\bigr)\int_{(0,t]}\frac{1}{1-\Delta
\Lambda(u)}\,dZ(u),
\end{eqnarray*}
which is a zero-mean Gaussian process with covariance function
\[
\bigl(1-F(s)\bigr)\bigl(1-F(t)\bigr)\int_{(0,s\wedge t]}\frac{1}{(1-\Delta
\Lambda(u))\overline{H}(u)}\,d\Lambda(u),
\]
then, by Theorem 5.2 of \citet{Arcones-LDP-04} and Theorem 3.1 of
\citet{Arcones-LDP-03}, we obtain the conclusion of
the theorem.
\end{pf}

\subsection{Moderate deviations for the empirical quantile processes}\label{sec4.3}

For a nondecreasing function $G\in D[a,b]$ and any $p\in\mathbb R$,
define $G^{-1}(p)=\inf\{x; G(x)\geq p\}.$ Let $D_1[a,b]$ denote the
set of all restrictions of distribution functions on $\mathbb R$ to
$[a,b]$ and let $D_2[a,b]$ denote the set of distribution functions
of measures that concentrate on $(a,b]$.

\begin{theorem}\label{quantile-processes-mdp-thm} Let $0<p<q<1$ be fixed and
let $F$ be a distribution function with continuous and positive
derivative $f$ on the interval
$[F^{-1}(p)-\varepsilon,F^{-1}(q)+\varepsilon]$ for some
$\varepsilon>0$. Let $F_n$ be the empirical distribution function
of an i.i.d. sample $X_1, \ldots, X_n$ of size $n$ from $F$. Then $
\{\frac{\sqrt{n}}{a(n)} (F_n^{-1}-F^{-1} ),n\geq1 \}
$ satisfies the LDP in $ l_\infty[p,q]$ with speed $a^2(n)$ and rate
function $I^{EQ}$ given by\looseness=-1
\[
I^{EQ}(\phi)=\inf \biggl\{I_{F}(\alpha); -\frac{\alpha
(F^{-1}(x))}{f(F^{-1}(x))} =\phi(x) \mbox{ for all } x\in[p,q]
\biggr\},
\]
where
\[
I_F(\alpha)=
\cases{\displaystyle \frac{1}{2}\int|\alpha'
_F(x)|^2\,dF(x), &\quad if $\alpha\ll F$ and
$\displaystyle\lim_{|t|\to\infty}|\alpha(t)|=0$, \vspace*{2pt}\cr\displaystyle
\infty, &\quad otherwise.
}
\]
\end{theorem}

\begin{pf} Applying Theorem 2 of \citet{Wu-94} to $L_n =\frac{1}{n}\sum_{i=1}^n \delta_{X_i}$,
\vspace*{-1pt} and
$\mathfrak{F} =\{(-\infty,x];x\in\mathbb R\}$, we know that $
\{\frac{\sqrt{n}}{a(n)} (F_n -F ),n\geq1 \} $ satisfies\vspace*{-1pt} the LDP on
$D(\mathbb R)$ with speed $ {a^2(n)} $ and rate function $I_F$. By
Lemma 3.9.23 of \citet{Vaar-Wellner-96}, it follows that
the inverse map $\Phi\dvtx G\mapsto G^{-1}$ as a map
$D_1[F^{-1}(p)-\varepsilon,F^{-1}(q)+\varepsilon]\mapsto
l_\infty[p,q]$ is Hadamard differentiable at $F$ tangentially to
$C[F^{-1}(p)-\varepsilon,F^{-1}(q)+\varepsilon]$, and the derivative
is the map $\alpha\mapsto-\alpha(F^{-1})/f(F^{-1}) $. Therefore,
by Theorem \ref{Delta-method-thm}, we conclude that $
\{\frac{\sqrt{n}}{a(n)} (F_n^{-1}-F^{-1} ),n\geq1 \} $ satisfies
the LDP in $ l_\infty[p,q]$ with speed $a^2(n)$ and the rate
function $I^{EQ}$.\vspace{-3pt}
\end{pf}

\subsection{Moderate deviations for the empirical copula processes}\label{sec4.4}

Let $BV_1^+(\mathbb R^2)$ denote the space of bivariate distribution
functions on $\mathbb R^2$. For $H\in BV_1^+(\mathbb R^2) $, set
$F(x)=H(x,\infty)$ and $G(y)=H(\infty,y)$.

Let $(X_1,Y_1),\ldots, (X_n,Y_n)$ be i.i.d. vectors with distribution
function $H$. The empirical estimator for the copula function
$C(u,v)=H(F^{-1}(u),G^{-1}(v))$ is defined by
$
C_n(u,v)=H_n(F_n^{-1}(u),G_n^{-1}(v)),
$
where $H_n$, $F_n$ and $G_n$ are the joint and marginal empirical
distributions of the observations.\vspace{-3pt}

\begin{theorem}\label{copula-mdp-thm} Let $0<p<q<1$ be fixed. Suppose that
$F$ and $G$ are
continuously differentiable on the intervals $
[F^{-1}(p)-\varepsilon,F^{-1}(q)+\varepsilon] $ and $
[G^{-1}(p)-\varepsilon,G^{-1}(q)+\varepsilon] $ with strictly
positive derivatives $f$ and $g$, respectively, for some
$\varepsilon>0$. Furthermore, assume that $\partial H/\partial x$
and $\partial H/\partial y$ exist and are continuous on the product
intervals. Then $ \{\frac{\sqrt{n}}{a(n)} (C_n-C ),n\geq1 \}
$
satisfies the LDP in $l_\infty([p,q]^2)$ with speed $a^2(n)$ and
rate function $I^C$ defined by\vspace{-1pt}
\[
I^{C}(\phi)=\inf \{I_{H}(\alpha); \Phi'_H(\alpha)=\phi\},\vspace{-3pt}
\]
where\vspace{-3pt}
\begin{eqnarray*}
\Phi'_H(\alpha)(u,v)&=&\alpha(F^{-1}(u),G^{-1}(v)) -\frac{\partial
H}{\partial
x}(F^{-1}(u),G^{-1}(v))\frac{\alpha(F^{-1}(u),\infty
)}{f(F^{-1}(u))}\\[-1pt]
&&{}-\frac{\partial H}{\partial
y}(F^{-1}(u),G^{-1}(v))\frac{\alpha(\infty,G^{-1}(u))}{g(G^{-1}(u))}.\vspace{-3pt}
\end{eqnarray*}
\end{theorem}

\begin{pf} By Theorem 2 of \citet{Wu-94}, we know
that\vspace{-1pt}
\[
P \Biggl(\frac{\sqrt{n}}{a(n)} \Biggl(\sum_{k=1}^n
\delta_{(X_k,Y_k)}\bigl((-\infty,x]\times
(-\infty,y]\bigr)-H(x,y) \Biggr)\in\cdot \Biggr)\vspace{-2pt}
\]
satisfies the LDP on $D(\mathbb R^2)$ with speed $ {a^2(n)} $ and
rate function defined as
\begin{eqnarray*}
I_H(\alpha)&=&\inf \biggl \{
\frac{1}{2}
\int\gamma^2(x,y)H(dx,dy);\
 \alpha(s,t)=\int\gamma(x,y) I_{\{x\leq s,y\leq t\}}H(dx,dy) \\[-3pt]
 &&\hspace*{138pt}\mbox{for each } (s,t)\in\mathbb R^2, \mbox{ and } \int\gamma
dH=0
 \biggr\}\\
&=&
\cases{\displaystyle \frac{1}{2}\int(\alpha'_H)^2(x,y)H(dx,dy), &\quad
 if $\alpha\ll H$ and
$\displaystyle\lim_{|s|,|t|\to\infty}|\alpha(s,t)|=0$, \vspace*{2pt}\cr\displaystyle
\infty, &\quad otherwise.
}
\end{eqnarray*}

Then, by Lemma 3.9.28 of \citet{Vaar-Wellner-96}, we
conclude that the map $\Phi\dvtx H\mapsto H(F^{-1},G^{-1})$ as a map
$BV_1^+(\mathbb R^2)\subset D(\overline{\mathbb R}^2)\mapsto
l_\infty([p,q]^2)$ is Hadamard differentiable at $H$ tangentially to
$C(\overline{\mathbb R}^2)$, and the derivative is $\Phi'_H$.
Therefore, it follows from Theorem \ref{Delta-method-thm} that,
$
\{\frac{\sqrt{n}}{a(n)} (C_n-C ),\break n\geq1 \}
$
satisfies the LDP in $l_\infty([p,q]^2)$ with speed $a^2(n)$ and
rate function~$I^C$ as defined in the theorem.
\end{pf}

\subsection{Moderate deviations for $M$-estimators}
$M$-estimators were first introduced by \citet{Huber-64}. Let $X$
be a
random variable taking its values in a measurable space
$(S,\mathcal S)$ with distribution $F$, let $X_{1}, \ldots, X_{n}$
be a random sample of $X$, and let $F_n$ denote the empirical
distribution function of~$X$. Let~$\Theta$ be a Borel subset of
$\mathbb R^d$. A $M$-estimator $\theta_n(X_1,\ldots, X_n)$ over the
function $g$ is a solution of
\[
\int g(x,\theta_n)\,dF_n(x)=\inf_{\theta\in\Theta}\int
g(x,\theta)\,dF_n(x).
\]
If $g(x,\theta)$ is differentiable with respect to $\theta$, then the
$M$-estimator $\theta_n(X_1,\ldots,\break X_n)$ may be defined as a solution
of the equation
\[
\int\nabla_\theta g(x,\theta_n)\,dF_n(x)=0,
\]
where $ \nabla_\theta g(x,\theta)=(\frac{\partial
g(x,\theta)}{\partial\theta^1},\ldots, \frac{\partial
g(x,\theta)}{\partial\theta^d})$. The detailed description on\vspace*{1pt}
$M$-estima\-tors can be found in \citet{Serfling-book-80}.

\citet{Jureckova-Kallenberg-Veraverbeke}, \citet
{Arcones-mdp-m-estimator-02} and \citet{Inglot-Kallenberg-03}
studied moderate
deviations for $M$-estimators. In this subsection, we study the
problem by the delta method. Let $\psi(x,\theta)=(\psi^1(x,\theta
), \ldots, \psi^d(x,\theta))\dvtx S\times
\Theta\mapsto\mathbb R^d$. We also need the following conditions.

(C1) $\psi(x,\theta)$ is continuous in $\theta$ for
each $x\in S$, and $\psi(x,\theta)$ is measurable in $x$ for each
$\theta\in\Theta$.

Define
\[
\Psi(\theta)=(\Psi^1(\theta),\ldots, \Psi^d(\theta))=E(\psi
(X,\theta))=\int\psi(x,\theta)\,dF(x)
, \qquad \theta\in\Theta,
\]
and
\[
\Psi_n(\theta)=(\Psi_n^1(\theta),\ldots, \Psi_n^d(\theta))=\frac
{1}{n}\sum_{i=1}^n \psi(X_i,\theta)=\int\psi(x,\theta
)\,dF_n(x), \qquad \theta\in\Theta.
\]

(C2)
$\Psi$ has a unique zero at $\theta_0$; there exists some $\eta>0$
such that $\bar{B}(\theta_0,\eta):=\{\theta\in\mathbb R^d;
 |\theta-\theta_0| \leq\eta\} \subset\Theta$ and $\Psi$ is
homeomorphism on $\bar{B}(\theta_0,\eta)$;
$\Psi$ is differentiable at $\theta_0$ with nonsingular derivative
$A\dvtx \mathbb R^d\mapsto\mathbb R^d$; and $E ( |\psi(X,\theta
)|^2 )<\infty$.

Let $C(\bar{B}(\theta_0,\eta))$ denote the space of continuous
$\mathbb R^d$-valued functions on $\bar{B}(\theta_0,\eta)$ and
define $\|f\|=\sup_{\theta\in\bar{B}(\theta_0,\eta)}|f(\theta)|$
for $f\in C(\bar{B}(\theta_0,\eta))$. Let $ \Psi_0(\theta)$ and
$\Psi_{0n}$ be the restrictions of $\Psi$ and $\Psi_n$ on
$\bar{B}(\theta_0,\eta)$, respectively.

(C3) $\{a(n),n\geq1\}$ satisfies
%
%
\begin{equation}\label{a(n)-increasing-c}
a(n )\nearrow\infty \quad \mbox{and} \quad \frac{a(n)}{\sqrt{n}} \searrow
0
\end{equation}
and $\{\psi(X_i,\theta),i\geq1\}$
satisfies
%
%
\begin{equation}\label{MD-LLN-c-M-est}
\frac{\sqrt{n}}{a(n)}\sup_{\theta\in\bar{B}(\theta_0,\eta
)}|\Psi_{n}(\theta)-\Psi(\theta)
|\stackrel{P}{\longrightarrow} 0
\end{equation}
and
%
%
\begin{equation}\label{Ledoux-A-W-c-M-est}
\limsup_{n\to\infty} \frac{1}{a^2(n)} \log \Bigl(nP \Bigl(
\sup_{\theta\in\bar{B}(\theta_0,\eta)}|\psi(X,\theta) |\geq
\sqrt{n}a(n) \Bigr) \Bigr)= -\infty.
\end{equation}

\begin{rmk}\label{Ledoux-A-W-c-examples}
Let $Y$ be a random variable taking its values in a Banach space and
$E(Y)=0$. If there exists a sequence of increasing
nonnegative functions $\{H_k, k\geq1 \}$ on $(0,+\infty)$
satisfying
%
%
\begin{equation}\label{h-c-1}
\lim_{u\to\infty}u^{-2}H_k(u)=+\infty,\qquad\lim_{k\to\infty
}\lim_{n\to\infty}\frac{1}{a^2(n)}\log
\frac{H_k(\sqrt{n}a(n))}{n}=+\infty,\hspace*{-40pt}
\end{equation}
and
%
%
\begin{equation}\label{h-c-2}
E(H_k(\|Y\|))<\infty\qquad\mbox{for any }
k\geq1,
\end{equation}
then
%
%
\begin{equation}\label{Ledoux-A-W-c} \limsup_{n\to\infty}
\frac{1}{a^2(n)} \log \bigl(nP \bigl (\|Y\|\geq\sqrt{n}a(n)
\bigr) \bigr)= -\infty.
\end{equation}
In particular [cf. \citet{chen91},
\citet{Ledoux-92}],
if for each $k\geq1$,
\[
E (\|Y\|^2({\log}
\|Y\|)^k )<+\infty,
\]
then (\ref{Ledoux-A-W-c})
holds for $a(n)=\sqrt{\log\log n}$;
if for each $k\geq1$,
\[
E(\|Y\|^k)<+\infty,
\]
then (\ref{Ledoux-A-W-c}) holds for
$a(n)=\sqrt{\log n}$;
if for some $1\leq p<2$, there exists some $\delta>0$
such that
%
%
\begin{equation}\label{integrable-c-1}
E (\exp \{\delta\|Y\|^{2-p} \} )<+\infty,
\end{equation}
then (\ref{Ledoux-A-W-c}) holds for $a(n)=o (n^{(2-p)/
{2p}} )$;
if for some $1< p<2$, and
%
%
\begin{equation}\label{integrable-c-2}
E (\exp \{\delta\|Y\|^{2-p} \} )<+\infty \qquad
\mbox{for all }
\delta>0,
\end{equation}
then (\ref{Ledoux-A-W-c}) holds for $a(n)=O (n^{(2-p)/
{2p}} )$.

In fact, by Chebychev's inequality,
\[
P\bigl(\|Y\| >\sqrt{n}a(n)\bigr)\leq\frac{E(H_k(\|Y\|))}{H_k(\sqrt{n}a(n))}.
\]
Hence, (\ref{h-c-1}) and (\ref{h-c-2}) yield (\ref{Ledoux-A-W-c}).
\end{rmk}

\begin{lem}[{[See Lemma 4.3 in \citet{Heesterman-Gill-92}]}]\label{M-est-H-dif-lem}
 Assume that \textup{(C1)} and
\textup{(C2)} hold.
Then there exists a neighborhood $V$ of $\Psi_0$ in $C(\bar{B}(\theta
_0,\eta))$
and a functional $\Phi\dvtx C(\bar{B}(\theta_0,\eta))\mapsto\bar
{B}(\theta_0,\eta)$ such that
\mbox{$f(\Phi(f))=0$}, for any $f\in V$,
and $\Phi$ is Hadamard differentiable at $\Psi_\eta$ with derivative
$
\Phi_{\Psi_0}'(f)=-A^{-1} f(\theta_0).
$
\end{lem}

\begin{theorem}\label{M-est-mdp-thm} Suppose that \textup{(C1)}, \textup{(C2)} and
\textup{(C3)} hold. Define
%
%
\begin{equation}\label{m-estimator-def}
\theta_n=\Phi(\Psi_{0n}).
\end{equation}
Then $\{\frac{\sqrt{n}}{a(n)}(\theta_n-\theta_0),n\geq1\}$
satisfies the LDP with speed
$a^2(n)$ and rate function
%
%
\begin{equation}\label{m-estimator-rate-f}
I^M(z)=\tfrac{1}{2} \langle Az, \Gamma^{-1} Az\rangle,
\end{equation}
where $\Gamma$ is the covariance of
$\psi(X,\theta_0)-\Phi(\theta_0)$, and
%
%
\begin{equation}\label{m-estimator-property-1}
\limsup_{n\to\infty}\frac{1}{a^2(n)}\log P \bigl(\Psi_{n}(\theta
_n)\not=0 \bigr)=-\infty.
\end{equation}
\end{theorem}

\begin{pf} Set $T=\{1,\ldots, d\}\times\bar{B}(\theta_0,\eta)$. Since
\[
T\times T\ni((i,s),(j,t))\mapsto d((i,s),(j,t)):=\bigl(\operatorname{Var}\bigl(\psi
^i(X_1,t)-\psi^j(X_1,s)\bigr)\bigr)^{1/2}
\]
is continuous on $ T\times T$ and $d((i,t),(i,t))=0$, then $(T,d)$
is totally bounded. Hence, under (C3), Theorem 2.8 in \citet
{Arcones-MDP-03} yields that\break
$\{\{\frac{\sqrt{n}}{a(n)}(\Psi_{0n}^i(\theta)-\Psi^i(\theta
)),(i,\theta)\in
T\},n\geq1\}$ satisfies the LDP in $l_\infty(T)$ with speed
$a^2(n)$ and rate function
\[
\hat{I}(f)= \tfrac{1}{2}\inf \bigl\{E(\alpha^2(X)); f(i,\theta
)=E\bigl(\alpha(X)\bigl(\psi^i(X,\theta)-\Psi^i(\theta)\bigr) \bigr)\bigr\}
\]
satisfying
\[
\limsup_{\lambda\to\infty}\frac{1}{\lambda}\inf \{\hat
{I}(f);\|f \|\geq\lambda \}=-\infty.
\]
Then, applying the classical contraction principle [see \citet
{Dembo-Z-98}, Theorem 4.2.1] to $l_\infty(\tilde{T})\ni f\to
(f(1,\theta),\ldots, f(d,\theta))\in l_\infty(\bar{B}(\theta
_0,\eta),\mathbb R^d)$, we obtain that
$\{\{\frac{\sqrt{n}}{a(n)}(\Psi_{0n}(\theta)-\Psi(\theta)),\theta
\in\bar{B}(\theta_0,\eta)\},n\geq1\}$ satisfies the LDP in $C(\bar
{B}(\theta_0,\eta))$ with speed $a^2(n)$ and rate
function
\[
I(f)=
\tfrac{1}{2}\inf \bigl\{E(\alpha^2(X)); f(\theta)=E\bigl(\alpha(X)\bigl(\psi
(X,\theta)-\Phi(\theta)\bigr)\bigr) \bigr\}.
\]
Therefore, we have
\[
\limsup_{n\to\infty}\frac{1}{a^2(n)}\log P (\Psi_{0n}\notin
V )=-\infty,
\]
and so (\ref{m-estimator-property-1}) holds. Then, by Theorem
\ref{Delta-method-thm}, we conclude that
$\{\frac{\sqrt{n}}{a(n)}(\theta_n-\theta),\break n\geq1\}$ satisfies the
LDP with speed $a^2(n)$ and rate function
\begin{eqnarray*}
I^M(z)&=&\tfrac{1}{2}\inf\bigl\{E(\alpha^2(X)),E\bigl(\alpha(X)\bigl(\psi(X,\theta
_0)-\Phi(\theta_0)\bigr)\bigr)=-Az\bigr\}\\
&=& \tfrac{1}{2} \langle Az, \Gamma^{-1} Az\rangle.
\end{eqnarray*}
\upqed
\end{pf}

\begin{rmk}
Comparing with Theorem 2.8 in \citet{Arcones-mdp-m-estimator-02},
in Theorem~\ref{M-est-mdp-thm}, we remove the condition
\[
\limsup_{n\to\infty}\frac{1}{a^2(n)}\log
P (|\theta_n-\theta_0|>\epsilon )=-\infty,
\]
which is required by \citet{Arcones-mdp-m-estimator-02}.
\end{rmk}

\subsection{Moderate deviations for $L$-statistics}\label{sec4.6}

Let $X_{1n}\leq X_{2n}\leq\cdots\leq X_{nn}$ be the order
statistics of a random sample $X_1, \ldots, X_n$ from a random
variable~$X$ with distribution function $F(x)$ and let $J$ be a
fixed score function on $(0,1)$. Also let $F_n$ be the empirical
distribution function of the sample. We consider the $L$-statistics
of the form
\[
L_n:=\sum_{i=1}^{n} X_{in} \int_{(i-1)/n}^{i/n} J(u)\,du=\int_0^{1}
F_n^{-1}(s)J(s)\,ds.
\]

\citet{GroeneboomOosterhoffRuymgaart} had obtained some large
deviations for $L$-statistics. The Cram\'{e}r type moderate deviations
for $L$-statistics had been studied in \citet{VV-MDP-L-stat-82},
\citet{Bentkus-Zitikis-90} and \citet{Aleskev-91}. In
this subsection, we study the moderate deviation principle for
$L$-statistics by the delta method.

Take $\XX=l_\infty(\mathbb R)$ and $\YY=\mathbb R$. Let $\DD_\Phi$
be the set of all distribution functions on $\mathbb R$, and set $
\DD_0=\{a(G-F);G\in\DD_\Phi,a\in\mathbb R\}$. Define $\Phi\dvtx
\DD_\Phi\mapsto\mathbb R$ as follows:
\[
\Phi(G) =\int_0^{1} G^{-1}(s)J(s)\,ds=\int_{-\infty}^{\infty}xJ(G(x))\,dG(x).
\]

Assume that $E(X^2)<\infty$. Set
$
m(J,F)= \int_{-\infty}^{\infty}xJ(F(x))\,dF(x),
$
and
\[
\sigma^2(J,F)=\int_{\mathbb R^2} J(F(x))J(F(y))\bigl(F(x\wedge
y)-F(x)F(y)\bigr)\,dx \,dy,
\]
where $x\wedge y=\min\{x,y\}$. We also assume $\sigma^2(J,F)>0$.

\begin{theorem}\label{mdp-L-stat-thm-1}
Suppose that the score function $J$ is trimmed near $0$ and $1$, that is,
$
J(u)=0, u\in[0,t_1)\cup(t_2,1]
$
where $0<t_1<t_2<1$. If $J$ is bounded and continuous a.e. Lebesgue
measure and a.e. $F^{-1}$, then\vspace*{-2pt}
$\{\frac{\sqrt{n}}{a(n)}(L_n-m(J,F)),\break n\geq1\}$ satisfies the LDP in
$ \mathbb R$ with speed $a^2(n)$ and rate function $
I^L(x)=\frac{x^2}{2 \sigma^2(J,F)}$.
\end{theorem}

\begin{pf} By Theorem 1 in \citet{Boos-79}, we have
\[
\lim_{\|G-F\|\to0}\frac{|\Phi(G)-\Phi(F)-\int
(F(x)-G(x))J(F(x))\,dx|}{\|G-F\|}=0.
\]
Therefore, for any $t_n\to0+$ and $H_n\to\alpha\in\DD_0$ with
$F+t_n H_n\in\DD_\Phi$,
\[
\lim_{n\to\infty} \biggl|\frac{|\Phi(F+t_n H_n)-\Phi(F)}{t_n}+
\int H_n(x)J(F(x))\,dx \biggr|=0,
\]
and so, $\Phi\dvtx \DD_\Phi\mapsto\mathbb R$ is Hadamard-differentiable
at $F$ tangentially to $\DD_0$ with respect to the uniform convergence,
and
$
\Phi_F'(\alpha)=- \int_{\mathbb R} \alpha(x) J(F(x))\,dx, \alpha
\in\DD_0.
$
By Theorem \ref{Delta-method-thm}, we conclude that
$\{\frac{\sqrt{n}}{a(n)}(L_n-m(J,F)),n\geq1\}$ satisfies the LDP in
$ \mathbb R$ with speed $a^2(n)$ and rate function $ I^L$ given by
\begin{eqnarray*}
I^L(y)=\inf \biggl\{I_F(\alpha); -\int_{\mathbb R} \alpha(x) J(F(x))\,dx=y \biggr\},
\end{eqnarray*}
which equals the rate function of\vspace*{-2pt}
\mbox{$\{-\frac{\sqrt{n}}{a(n)}\int_{\mathbb R}(F_n(x)-F(x)) J(F(x))\,dx,n\geq1\}$}, that is,
$
I^L(y)=\frac{y^2}{2 \sigma^2(J,F)}.
$
\end{pf}

Now, let us remove the trimming restrictions on $J$.
Set
\begin{eqnarray*}
\tilde{\DD}_\Phi&=& \biggl \{
 \tilde{G}(x)=G(x)I_{(-\infty
,0)}(x)+\bigl(G(x)-1\bigr)I_{[0,\infty)}(x);\\
&&\hspace*{100pt} G\in\DD_\Phi,\int|x|\,dG(x)<\infty
 \biggr\}
\end{eqnarray*}
and
$
\tilde{\DD}_0= \{a( \tilde{G} -\tilde{F} )\equiv a(G-F); a\in
\mathbb R,\tilde{G}\in
\tilde{\DD}_\Phi \}.
$
Then $\tilde{\DD}_\Phi,\tilde{\DD}_0\subset
L^1(\mathbb R)$. Define $\tilde{\Phi}\dvtx \tilde{\DD}_\Phi\mapsto
\mathbb R$ by
$\tilde{\Phi}(\tilde{G})=\Phi(G)$ for all $\tilde{G}\in\tilde
{\DD}_\Phi$.

\begin{lem}\label{L-stat-H-diff-lemmma} If $J$ is Lipschitz continuous on
$[0,1]$, then\vspace*{1pt} $\tilde{\Phi}\dvtx \tilde{\DD}_\Phi\mapsto\mathbb R$
is Ha\-damard-differentiable at $\tilde{F}$ tangentially to $\tilde{\DD
}_0$ with
respect to $L^1$-convergen\-ce, and
\[
\tilde{\Phi}_{\tilde{F}}'(\alpha)=- \int_{\mathbb R} \alpha(x) J(F(x))\,dx, \qquad \alpha\in\tilde{\DD}_0.
\]
\end{lem}
\vfill\eject

\begin{pf} By integration by parts, we can write [cf. \citet{Boos-79},
\citet{ShaoJun-89}]
\[
\tilde{\Phi}(\tilde{G}) -\tilde{\Phi}(\tilde{F}) +\int_{\mathbb
R} \bigl(G(x)-F(x)\bigr)J(F(x))\,dx=R(G,F) \qquad \mbox{for any }\tilde{G}\in\tilde
{\DD}_\Phi,
\]
where
$
R(G,F)=\int_{\mathbb R} W_{G,F}(x)(G(x)-F(x))\,dx,
$
and
\[
W_{G,F}(x)=
\cases{\displaystyle \frac{\int_{F(x)}^{G(x)}(
J(t)-J(F(x)))\,dt}{G(x)-F(x)},&\quad if $G(x)\not=F(x)$,\cr\displaystyle
0,&\quad if $G(x)=F(x)$.
}
\]
By the Lipschitz continuity of $J$, there exists a constant $C>0$
such that
\[
|R(G,F)|\leq C \int_{\mathbb R} \bigl( G (x)- F (x)\bigr)^2\,dx=C \int_{\mathbb
R} \bigl(\tilde{G}(x)-\tilde{F}(x)\bigr)^2\,dx.
\]

For any $t_n\to0+$ and $H_n\to\alpha\in\tilde{\DD}_0$
in $(L^1(\mathbb R),\|\cdot\|_{L^1})$ with $\tilde{F}+t_n H_n\in
\tilde{\DD}_\Phi$, then $|H_n|\leq2/t_n$ and
\[
\int_{\mathbb R}|H_n(x)-\alpha(x)|^2\,dx \leq(\|\alpha\|+2/t_n)\int
_{\mathbb R}|H_n(x)-\alpha(x)|\,dx,
\]
where $\|\alpha\|=\sup_{x\in\mathbb R} |\alpha(x)|$.
Therefore,
\begin{eqnarray*}
&&\frac{1}{t_n}\int_{\mathbb R} \bigl(\tilde{F}(x)+t_n H_n(x)-\tilde
{F}(x)\bigr)^2\,dx\\
&& \qquad \leq 2 t_n\int_{\mathbb R}|H_n(x)-\alpha(x)|^2\,dx+2t_n\int_{\mathbb
R}| \alpha(x)|^2\,dx\to0,
\end{eqnarray*}
and so
\[
\lim_{n\to\infty} \biggl|\frac{\tilde{\Phi}(\tilde{F}+t_n
H_n)-\tilde{\Phi}(\tilde{F})}{t_n}+ \int
H_n(x)J(F(x))\,dx \biggr|=0,
\]
which yields that $\tilde{\Phi} $ is
Hadamard-differentiable at $\tilde{F}$ tangentially to $\tilde{\DD
}_0$ with
respect to $L^1$-convergence, and
$
\tilde{\Phi}_{\tilde{F}}'(\alpha)=- \int_{\mathbb R} \alpha(x) J(F(x))\,dx $.
\end{pf}

\begin{lem}\label{Totally-b-lem}
Let $X$ be a random variable with
values in a separable Banach space $\mathbf{B}$ and $E(\|X\|^2)<\infty$. Then
$(\mathbf{B}_1^*,d)$ is totally bounded,
where $\mathbf{B}_1^*$ is the unit ball of the dual space $\mathbf{B}^*$ of
$\mathbf{B}$, and
\[
d(g,h)= \bigl(E \bigl( \bigl(g\bigl(X-E(X)\bigr)-h\bigl(X-E(X)\bigr) \bigr)^2
\bigr) \bigr)^{1/2}, \qquad g,h\in\mathbf{B}_1^*.
\]
\end{lem}
\begin{pf}
Noting $|g(X-E(X))-h(X-E(X))|\leq2\|X-E(X)\|$ for all $g,h\in\mathbf{B}_1^*$ and $E(\|X-E(X)\|^2)<\infty$, by the dominated convergence
theorem, we know that the function $(g,h)\mapsto d(g,h)$ is
continuous from $\mathbf{B}_1^*\times\mathbf{B}_1^*$ to $\mathbb R$ with
respect to $w^*$-topology. Let $ d^*$ denote a compatible metric on
$(B_1^*,w^*)$. Since $\mathbf{B}_1^*$ is $w^*$-compact and $d(g,g)=0$,
then, for any $\epsilon>0$, there exists some $\delta>0$ such that
$
d(g,h)<\epsilon, \mbox{ if } d^*(g,h)<\delta.
$
Choose finite points $h_1,\ldots, h_m\in\mathbf{B}_1^*$ such that
$B_1^*\subset\bigcup_{i=1}^m \{g;d^*(g,h_i)<\delta\}$, then
$B_1^*\subset\bigcup_{i=1}^m \{g;d(g,h_i)<\epsilon\}$. Therefore,
$(\mathbf{B}_1^*,d)$ is totally bounded.
\end{pf}

Define
\[
\Lambda_{2,1}(X)=\int_0^\infty\sqrt{P(|X|>t)}\,dt.
\]
Then [cf. \citet{Barrio-Gine-Matran-99}, page~1014],
$\Lambda_{2,1}(X)<\infty$ if and only if
$
\int_{-\infty}^\infty\sqrt{F(x)(1-F(x))}\,dx<\infty.
$

\begin{lem}\label{mdp-W-metric} Assume that $\Lambda_{2,1}(X)<\infty
$. If (\ref{a(n)-increasing-c}) holds and
%
%
\begin{equation}\label{Ledoux-A-W-c-L}
\limsup_{n\to\infty} \frac{1}{a^2(n)} \log \bigl(nP \bigl (
|X|\geq
\sqrt{n}a(n) \bigr) \bigr)= -\infty,
\end{equation}
then $\{\frac{\sqrt{n}}{a(n)}(F_n-F)=\frac{\sqrt{n}}{a(n)}(\tilde
{F}_n-\tilde{F}),n\geq1\}\subset\tilde{\DD}_0$
satisfies the LDP in $(L^1(\mathbb R),\break \|\cdot\|_{L^1})$ with speed $a^2(n)$
and rate function $I_F$.
\end{lem}
\begin{pf}
Set $\xi_i=I_{\{X_i\leq x\}}-F(x),x\in\mathbb R$, then
\[
\|\xi_i\|_{L^1}=2 \biggl(X_iF(X_i)-\int_{(-\infty,X_i)}x \,dF(x) \biggr).
\]
Therefore, the condition of the lemma implies
\[
\limsup_{n\to\infty} \frac{1}{a^2(n)} \log \bigl(nP \bigl(
\|\xi_1\|_{L^1}\geq\sqrt{n}a(n) \bigr) \bigr)= -\infty,
\]
and by Theorem 2.1(b) of \citet{Barrio-Gine-Matran-99},
we also\vspace*{-2pt} have
$
\frac{1}{a(n)}\|\sum_{i=1}^n\xi_i\|_{L^1}\stackrel
{P}\longrightarrow0.
$
By Lemma \ref{Totally-b-lem},
$(\mathbf{B}_1^*,d)$ is totally bounded, where
\[
\mathbf{B}_1^*:=\Bigl\{g\in L^\infty; \|g\|_\infty:=\operatorname{esssup}\limits_{x\in
\mathbb R}|g(x)|\leq1\Bigr\}
\]
and
\[
d(g,h)= \biggl(E \biggl( \biggl(\int_{\mathbb R}\bigl(g(x)-h(x)\bigr)\xi
_1(x)\,dx \biggr)^2 \biggr) \biggr)^{1/2}.
\]
Therefore, by Theorem 2.8 in \citet{Arcones-MDP-03}, the
conclusion of the
lemma holds.
\end{pf}

By Lemmas \ref{mdp-W-metric} and \ref{L-stat-H-diff-lemmma} and
Theorem \ref{Delta-method-thm}, we obtain the following result.

\begin{theorem}\label{mdp-L-stat-thm-2}
Assume that $\Lambda_{2,1}(X)<\infty$, (\ref{a(n)-increasing-c}) and
(\ref{Ledoux-A-W-c-L}) hold.
If $J$ is Lipschitz continuous on $[0,1]$, then\vspace*{-2pt} $\{\frac{\sqrt
{n}}{a(n)}(L_n-m(J,F)),n\geq1\}$ satisfies the LDP in $ \mathbb R$
with speed $a^2(n)$ and rate
function $
I^L(x)=\frac{x^2}{2 \sigma^2(J,F)}.
$
\end{theorem}

\begin{rmk} From Remark \ref{Ledoux-A-W-c-examples}, the moment
condition in Theorem~\ref{mdp-L-stat-thm-2} is weaker than the
conditions given in \citet{VV-MDP-L-stat-82},
\citet{Bentkus-Zitikis-90} and \citet{Aleskev-91}. In
particular, if\break
$E(|X|^{2+\delta})<\infty$ and $a(n)=\sqrt{\log\log n}$, then the
condition of Lemma \ref{mdp-W-metric} is valid, and so, for any $r>0$,
\[
\lim_{n\to\infty} \frac{1}{\log\log n} \log P \Biggl(\sqrt{\frac
{n}{\log\log n}}|L_n-m(J,F)|\geq r \Biggr)=-\frac{r^2}{2\sigma^2(J,F)}.
\]
\end{rmk}

\section{Application: Statistical hypothesis testing}\label{sec5}

In this section, we applied the moderate deviations to hypothesis
testing problems. We only consider the right-censored data model.
The method can be applied to other models.

Let $F$ be the unknown distribution function in the right-censored
data model considered in Section~\ref{sec4.2} and let $\hat{F}_n$ be the
Kaplan--Meier estimator of $F$. Consider the following hypothesis
testing:
\[
H_0\dvtx F=F_0\quad\mbox{and} \quad H_1\dvtx F=F_1,
\]
where $F_0$ and $F_1$ are two distribution functions such that
$F_0(x_0)\not=F_1(x_0)$ for some $x_0\in[0, \tau]$. Similar to the
Kolmogorov--Smirnov test, we take the Kaplan--Meier statistic
$T_n:=\sup_{x\in[0,\tau]}|\hat{F}_n(x)-F_0(x)|$ as test statistic.
Suppose that the rejection region for testing the null hypothesis
$H_0$ against $H_1$ is $\{\frac{\sqrt{n}}{a(n)}T_n\geq c\}$,\vspace*{1pt} where
$c$ is a positive constant. Then the probability $\alpha_n$ of Type
I error and the probability $\beta_n$ of Type II error are
\[
\alpha_n=P \biggl(\frac{\sqrt{n}}{a(n)}T_n\geq
c \Big|F=F_0 \biggr) \quad \mbox{and} \quad \beta_n=P \biggl(\frac{\sqrt{n}}{a(n)}T_n<
c \Big|F=F_1 \biggr),
\]
respectively. Then
\begin{eqnarray*}
\beta_n &\leq& P \biggl(\frac{\sqrt{n}}{a(n)}\sup_{x\in[0,\tau]}|\hat
{F}_n(x)-F_1(x)|\\
&&\hphantom{P \biggl(}\geq\frac{\sqrt{n}}{a(n)}\sup_{x\in[0,\tau
]}|F_0(x)-F_1(x)|-c \Big|F=F_1 \biggr).
\end{eqnarray*}
Therefore, Theorem \ref{Kaplan-Meier-mdp-thm} implies that
\[
\lim_{n\to\infty}\frac{1}{a^2(n)}\log
\alpha_n=- \frac{c^2}{2\sigma_{KM}^2},
\qquad
\lim_{n\to\infty}\frac{1}{a^2(n)}\log
\beta_n=- \infty,
\]
where
\begin{eqnarray*}
\sigma_{KM}^2&=&\sup_{t\in[0,\tau]} \bigl(1-F_0(t)\bigr)^2\int_{(0,
t]}\frac{1}{(1-\Delta\Lambda(u))\overline{H}_0(u)}\,d\Lambda_0(u),
\\
\Lambda_0(t)&=&\int_{[0,t]}\frac{1}{1-F_0(s-)}\,dF_0(s), \qquad
\overline{H}_0(t)=P(Z \geq t|F=F_0),
\end{eqnarray*}
and $Z$ is as defined in Section~\ref{sec4.2}.

The above result tells us that if the rejection region for the test
is $\{\frac{\sqrt{n}}{a(n)}T_n\geq c\}$, then the probability of
Type I error tends to $0$ with decay speed
\[
\exp\{- c^2
a^2(n)/(2\sigma_{KM}^2) \},
\]
and the probability of Type II error
tends to $0$ with decay speed $\exp\{-r a^2(n)\}$ for all $r>0$.

\section{Concluding remarks}\label{sec6}
This article discussed the large
deviations of transformed statistics. For
the problem, an extended contraction principle was developed and a
general delta method in large deviation theory was proposed. The
new method was used to establish the moderate deviation principles
for the Wilcoxon statistic, the Kaplan--Meier estimator, the
empirical quantile estimator and the empirical copula estimator,
which have not been addressed in the literature. The proposed method
was also used to improve the existing moderate deviation results for
$M$-estimators and $L$-statistics, where our proofs are different
from others but simpler by the new method. Moreover, our moderate
deviation results are very useful for statistical hypothesis
testing. As shown in Section~\ref{sec5}, a moderate deviation result can be
used to construct a test of a statistical hypothesis such that the
probabilities of both Type I and Type II errors tend to $0$ with an
exponentially decay speed as $n\to\infty$.

Note that the asymptotics for multivariate trimming and general $Z$-estima\-tors have been studied by using the delta method in a weak
convergence; see \citet{Nolan-92} and \citet{Vaar-Wellner-96}.
Similar to those presented in Section~\ref{sec4}, the moderate deviations for
these estimators can be established by using the proposed delta
method in large deviations.

These applications show that the proposed method is very powerful
for deriving moderate deviation principles on estimators. The method
will play an important role in large sample theory of statistics
like the functional delta method in weak convergence. Theoretically
speaking, we can apply the proposed delta method to obtain moderate
deviations for estimators where the classical delta method can be
applied.

%
\begin{appendix}

\section*{Appendix: Proof of the extended contraction principle}\label{appm}

\textit{Step 1}. First of all, let us prove $\{I<\infty\}\subset
\DD_\infty$, where $\DD_{\infty}$ denotes the set of all $x$ for
which there exists a sequence $x_{n}$ with $x_{n}\in\DD_{n}$ and
$x_{n}\to x$.

In fact, by the definition of $\DD_\infty$, $x\in\DD_\infty$ if and
only if for any $k\geq1$, there exists a positive integer $n_k$
such that $B_d(x,1/k)\cap\DD_n \not=\varnothing$ for all $n\geq n_k$,
where $B_d(x,1/k)=\{y\in\XX;d(y,x)<1/k\}$. Therefore, for any
$x\notin\DD_\infty$, there exist an open neighborhood $U$ of $x$ and
a subsequence $\{\DD_{n_k},k\geq1\}$ such that $\DD_{n_k}\cap
U=\varnothing$ for all $k\geq1$. Then by the lower bound of the
large deviations for $\{X_n,n\geq1\}$, we have
\[
-\infty=\liminf_{k\to\infty}\frac{1}{\lambda(n_k)}\log
{P_{n_k}}_*(X_{n_k}\in U)\geq-I(x),
\]
which implies $\{I<\infty\}\subset\DD_\infty$, where ${P_{n_k}}_*$
is the inner measure correspon\-ding to $P_{n_k}$ as defined in
Section~\ref{sec2}.

 \textit{Step 2}. Let us prove that
if some subsequence $x_{n_k}\to x\in\{I<\infty\}$ with $x_{n_k}\in
\DD_{n_k}$, then $f_{n_k}(x_{n_k})\to f(x)$ and the restriction of the
function $f$ to $\{I<\infty\} $ is continuous.

The proof is
similar to that of the extended mapping theorem [see Theorem~1.11.1
in \citet{Vaar-Wellner-96}], which is given below.
Let a subsequence $x_{n_k}\to x\in\{I<\infty\}$ be
given. Since $x\in\DD_\infty$, there exists a~sequence $y_n\to x$
with $y_n\in\DD_n$ for each $n\geq1$. Define
$
x_n=x_{n}I_{\{n_k,k\geq1\}}(n)+ y_nI_{\{n_k,k\geq1\}^c}(n)$.
Then $x_n\in\DD_n$ for each $n\geq1$ and $x_n\to x$. Therefore, by
 condition (ii), $f_n(x_n)\to f(x)$, and so $f_{n_k}(x_{n_k})\to
f(x)$. To prove the continuity of $f$ on $\{I<\infty\}$,
let $x_m\to x$ in $ \{I<\infty\}$. For every $m$, there is
a~sequence $x_{m,n}\in\DD_{n}$ with $x_{m,n}\to x_m$ as $n\to\infty$.
Since $x_m\in\{I<\infty\}$, then $f_{n}(x_{m,n})\to
f(x_m)$ as $n\to\infty$. For every~$m$, take $n_m$ such that
$n_{m}$ is increasing with $m$ satisfying
$
d(x_{m,n_{m}},x_m)<1/m $ and $\rho(f_{n_{m}}(x_{m,n_{m}}), f(x_m))<1/m.
$
Then $x_{m,n_{m}}\to x$, and by the first conclusion in Step 2,
$f_{n_{m}}(x_{m,n_{m}})\to f(x)$ as $m\to\infty$. This yields
$f(x_m)\to f(x)$.

 \textit{Step 3}. Let us prove that $[I_{ f}\le L] = f([I\le L])$
for any $ L\geq0$ and $I_f$ is inf-compact, that is, for any
$L\in[0,+\infty)$, $[I_{ f}\le L]$ is compact. This can be shown by
the continuity of $f|_{\{I<\infty\}}$ obtained in Step 2.

 \textit{Step 4}. Next, we show the upper bound of large
deviations.

Let $F$ be a closed subset in $\YY$. Then, using the
arguments similar to the proof of the extended continuous mapping
theorem [see Theorem 1.11.1 in \citet{Vaar-Wellner-96}],
%
%
\begin{equation}\label{ext-contraction-P-thm-ULD-eq-4}
\bigcap_{n=1}^\infty\overline{\bigcup_{m=n}^\infty f_m^{-1}(F)}\subset
f^{-1}(F)\cup( \{I<\infty\})^c.
\end{equation}

 Now for every fixed $k$, by the large deviation principle of
$\{X_n,n\geq1\}$, for each $L>0$, there exists a compact subset $K_L$
such that for any $\delta>0$,
\[
\limsup_{n\to\infty}\frac{1}{\lambda(n)}\log P_n^* \bigl(X_n\in
(K_L^\delta)^c \bigr)\leq-L,
\]
and so
\begin{eqnarray*}
&&\limsup_{n\to\infty}\frac{1}{\lambda(n)}\log P_n^*\bigl(f_n(X_n)\in
F\bigr)\\
&& \qquad \leq\limsup_{n\to\infty}\frac{1}{\lambda(n)}\log P_n ^*
\Biggl(X_n\in
\overline{\bigcup_{m=k}^\infty f_m^{-1}(F)} \Biggr)\\
&& \qquad \leq\max\biggl \{ -\inf_{x\in
\overline{K_L^\delta}\cap\overline{\bigcup_{m=k}^\infty f_m^{-1}(F)}}
I(x),-L \biggr\},
\end{eqnarray*}
where $ K_L^\delta=\{y; d(y,x)< \delta \mbox{ for some } x\in K_L
\}$ and $P_n^*$ is the outer measure corresponding to $P_n$ as
defined in Section~\ref{sec2}. Since $K_L $ is compact and $I$ is lower
semi-continuous, then, when $\delta\downarrow0$,
\[
\inf_{x\in\overline{K_L^\delta}\cap\overline{\bigcup_{m=k}^\infty
f_m^{-1}(F)}} I(x)\uparrow\inf_{x\in
K_L\cap\overline{\bigcup_{m=k}^\infty f_m^{-1}(F)}} I(x).
\]
Hence it follows that
\[
\limsup_{n\to\infty}\frac{1}{\lambda(n)}\log P_n ^*\bigl(f_n(X_n)\in F\bigr)
\leq\max \biggl\{ -\inf_{x\in
K_L\cap\overline{\bigcup_{m=k}^\infty f_m^{-1}(F)}} I(x),-L \biggr\}.
\]
Choose a sequence $x_k\in K_L\cap\overline{\bigcup_{m=k}^\infty
f_m^{-1}(F)}$, $k\geq1$ such that
$
I(x_k)=\break \inf_{x\in K_L\cap\overline{\bigcup_{m=k}^\infty f_m^{-1}(F)}}
I(x),
$
and then choose a subsequence $\{x_{k_m},m\geq1\}$ and $x_0\in K_L$
such that $x_{k_m} \to x_0$. Then we have
\[
x_0\in K_L\cap\Biggl(\bigcap_{k=1}^\infty\overline{\bigcup_{m=k}^\infty
f_m^{-1}(F)}\Biggr)\subset K_L\cap\bigl(f^{-1}(F)\cup(\{I<\infty\})\bigr).
\]
Letting
$k\to\infty$, we have
\[
\liminf_{k\to\infty}I(x_k)\geq I(x_0) \geq\inf_{x\in
K_L\cap(f^{-1}(F)\cup(\{I<\infty\}) ^c)} I(x)\geq
\inf_{x\in f^{-1}(F) }I(x).
\]
Now letting $L \to\infty$, we conclude that
\[
\limsup_{n\to\infty}\frac{1}{\lambda(n)}\log P_n ^*\bigl(f_n(X_n)\in
F\bigr) \leq
-\inf_{x\in
f ^{-1}(F) } I(x)=-\inf_{x\in F} I_f(x).
\]

 \textit{Step 5}. Finally, we show the lower bound of large
deviations: for any $y_0\in\YY$ with $I_f(y_0)<\infty$,
\[
\liminf_{n\to\infty}\frac{1}{\lambda(n)}\log{P_n}_*\bigl(f_n(X_n)\in
B(y_0,\delta)\bigr)\geq- I_f(y_0).
\]

For any $a>I_f(y_0)$, there is some $x_0\in\XX$ with $f(x_0)=y_0$ and
$I(x_0)<a$. For any $\delta>0$, set
$
B(\delta)=B_\rho(y_0,\delta)=\{y\in\YY;\rho(y_0,y)<\delta\}
$
and
$F(\delta)=B(\delta)^c$. Then, by
(\ref{ext-contraction-P-thm-ULD-eq-4}), we
have
%
%
\begin{equation}\label{ext-contraction-P-thm-LLD-eq-1}
\bigcup_{n=1}^\infty\overline{\bigcup_{m=n}^\infty
f_m^{-1}(F(\delta))}^{\,c}\supset f^{-1}(B(\delta))\cap
(\{I<\infty\})\ni x_0.
\end{equation}
Now for every fixed $k$, by the large deviation principle of
$\{X_n\}$, we have
\begin{eqnarray*}
&&\liminf_{n\to\infty}\frac{1}{\lambda(n)}\log{P_n}_*
\bigl(f_n(X_n)\in
B(\delta)\bigr)\\
&& \quad \geq
\liminf_{n\to\infty}\frac{1}{\lambda(n)}\log{P_n}_* \Biggl(X_n\in
\overline{\bigcup_{m=k}^\infty f_m^{-1}(F(\delta))}^{\,c} \Biggr)\geq
-\inf_{x\in\overline{\bigcup_{m=k}^\infty
f_m^{-1}(F(\delta))}^{\,c}} I(x).
\end{eqnarray*}
Since
$
x_0\in f^{-1}(B(\delta)) \subset
\bigcup_{n=1}^\infty\overline{\bigcup_{m=n}^\infty f_m^{-1}(F(\delta))}^{\,c}
,
$
there is some $k\geq1$ such that $x_0\in
\overline{\bigcup_{m=k}^\infty f_m^{-1}(F(\delta))}^{\,c}$. Therefore,
\[
\liminf_{n\to\infty}\frac{1}{\lambda(n)}\log{P_n}_* \bigl(f_n(X_n)\in
B(\delta)\bigr)\geq- I(x_0)>-a.
\]
Letting $a\downarrow I_f(y_0)$, we obtain the lower bound of large
deviations.

\begin{rmk}
When $\DD_n=\XX$ for all $n\geq1$, the continuity of $f$ can be
proved directly by the following property [see Theorem 2.1 in
\citet{Arcones-MDP-03}]:
Given $\epsilon>0$, for any $x_0\in\{I<\infty\}$, there are
$\delta>0$ and a positive integer $n_0$ such that for all $n\geq
n_0$,
$
f_n(B(x_0,\delta))\subset B(f(x_0),\epsilon).
$
However, when $\DD_n\not=\XX$, $f_n(B(x_0,\delta))$ is not well
defined since $B(x_0,\delta)\not\subset\DD_n$. Thus, the above
property cannot be used for proving the continuity of $f$ in this
case.
\end{rmk}
\end{appendix}

\section*{Acknowledgments}
The authors are very grateful to the
Editors, Professors Tony Cai and Bernard W. Silverman, the Associate
Editor and the three referees for their many valuable comments and
suggestions that greatly improved the paper. The authors also wish
to thank Professor Qiman Shao and Professor Liming Wu for their
useful discussions and suggestions.


%

\printaddresses

\end{document}